\documentclass[psamsfonts]{amsart}
\usepackage{amssymb}

\newcommand{\jj}{\vee}
\newcommand{\mm}{\wedge}
\newcommand{\JJ}{\bigvee}
\newcommand{\MM}{\bigwedge}
\newcommand{\JJm}[2]{\JJ(\,#1\mid#2\,)}
\newcommand{\MMm}[2]{\MM(\,#1\mid#2\,)}

\newcommand{\uu}{\cup}
\newcommand{\ii}{\cap}
\newcommand{\UU}{\bigcup}
\newcommand{\II}{\bigcap}
\newcommand{\UUm}[2]{\UU(\,#1\mid#2\,)}
\newcommand{\IIm}[2]{\II(\,#1\mid#2\,)}

\newcommand{\ci}{\subseteq}
\newcommand{\sci}{\subset}
\newcommand{\ce}{\supseteq}
\newcommand{\nin}{\notin}
\newcommand{\set}[1]{\{#1\}}
\newcommand{\setm}[2]{\{\,#1\mid#2\,\}}
\def\vv<#1>{\langle#1\rangle}

\newcommand{\ga}{\alpha}
\newcommand{\gb}{\beta}
\renewcommand{\gg}{\gamma}
\newcommand{\gd}{\delta}
\renewcommand{\ge}{\varepsilon}
\newcommand{\gi}{\iota}
\newcommand{\go}{\omega}
\newcommand{\gw}{\omega}
\newcommand{\gF}{\Phi}
\newcommand{\gQ}{\Theta}

\newcommand{\tbf}{\textbf}
\newcommand{\tup}{\textup}
\newcommand{\mrm}{\mathrm}
\newcommand{\E}[1]{\mathcal{#1}}

\newcommand{\ol}[1]{\overline{#1}}
\newcommand{\ul}[1]{\underline{#1}}
\newcommand{\q}{\quad}
\newcommand{\qq}{\qquad}
\newcommand{\iso}{\cong}

\theoremstyle{plain}
\newtheorem{theorem}{Theorem}
\newtheorem{lemma}{Lemma}[section]
\newtheorem{proposition}[lemma]{Proposition}
\newtheorem{corollary}[lemma]{Corollary}
\newtheorem{claim}{Claim}
\newtheorem*{stat}{\name}
\newcommand{\name}{testing}
\theoremstyle{definition}
\newtheorem{definition}[lemma]{Definition}
\newtheorem{notation}[lemma]{Notation}
\newtheorem{example}[lemma]{Example}
\newtheorem{problem}{Problem}
\newtheorem*{note}{Note}
\newtheorem{remark}[lemma]{Remark}

\newcommand{\congtimes}{\mathbin{\square}}
\newcommand{\tim}[1]{\mathbin{\congtimes_{#1}}}
\newcommand{\bboxtimes}{\mathbin{\ol{\congtimes}}}
\newcommand{\congtens}{\mathbin{\odot}}
\newcommand{\odc}[1]{\mathbin{\congtens_{#1}}}
\newcommand{\ootimes}{\mathbin{\ol{\otimes}}}
\newcommand{\jz}{$\set{\jj, 0}$}
\newcommand{\dd}{{\mathrm{d}}}

\newenvironment{all}[1]{\renewcommand{\name}{#1}\begin{stat}}
                        {\end{stat}}

\DeclareMathOperator{\Con}{Con}
\DeclareMathOperator{\Pow}{Pow}
\DeclareMathOperator{\Conc}{Con_c}
\DeclareMathOperator{\ConL}{Con^L}
\DeclareMathOperator{\Id}{Id}

\begin{document}

\title[Tensor products of semilattices]
{Tensor products of semilattices with zero, revisited}

\author{G.~Gr\"atzer}
\thanks{The research of the first author was
        supported by the
        NSERC of Canada.}
\address{Department of Mathematics\\
	  University of Manitoba\\
	  Winnipeg MN, R3T 2N2\\
	  Canada}
 \email{gratzer@cc.umanitoba.ca}
 \urladdr{http://www.maths.umanitoba.ca/homepages/gratzer.html/}

\author{F.~Wehrung}
\address{C.N.R.S., E.S.A. 6081\\
   D\'epartement de Math\'ematiques\\
	  Universit\'e de Caen\\
	  14032 Caen Cedex\\
   France} \email{wehrung@math.unicaen.fr}
\urladdr{http://www.math.unicaen.fr/\~wehrung}

\date{January 26, 1998}

 \begin{abstract}
 Let $A$ and $B$ be lattices with zero. The classical tensor product,
$A\otimes B$, of $A$ and $B$ as join-semilattices with zero is a
join-semilattice with zero; it is, in general, not a lattice.  We define a
very natural condition: $A\otimes B$ is \emph{capped} (that is, every element
is a finite union of pure tensors) under which the tensor product is always a
lattice.

 Let $\Conc L$ denote the join-semilattice with zero of compact congruences
of a lattice $L$. Our main result is that the following isomorphism holds for
any capped tensor product:
 \[
   \Conc A \otimes \Conc B \iso \Conc (A \otimes B).
 \]
 This generalizes from finite lattices to arbitrary lattices the main
result of a joint paper by the first author, H. Lakser, and R. W. Quackenbush.
 \end{abstract}

\maketitle

\section{Introduction}\label{S:Introduction}
 The construction of tensor products of modules over a commutative ring has an
obvious analogue for join-semilattices with zero. This construction was
introduced in J.~Anderson and N.~Kimura \cite{AK78}, G.~A. Fraser
\cite{Fras76},
and Z.~Shmuley \cite{zS74}.  If $A$ and~$B$ are semilattices with zero, we
denote by $A \otimes B$ the tensor product of $A$ and~$B$.

 While the tensor product is defined for semilattices with zero, it becomes,
somehow mysteriously, really interesting for \emph{lattices}. In many
cases, the
tensor product of two lattices with zero is a lattice, for example, if both
lattices are finite, see J.~Anderson and N.~Kimura \cite{AK78}.

 The deepest result in this field was obtained in G. Gr\"atzer, H. Lakser, and
R.~W. Quackenbush \cite{GLQu81}.  This paper was motivated by the paper of
E. T.
Schmidt \cite{tS68} in which it is proved that the congruence lattice of
$M_3[D]$, where $D$ is a bounded distributive lattice, is isomorphic to the
congruence lattice of $D$.  Since $M_3[D]$ can be viewed as a tensor product of
$M_3$ and $D$, the following result is a far reaching and surprising
generalization of Schmidt's result:

 \begin{all}{Main Result of \cite{GLQu81}}
 Let $A$ and~$B$ be finite lattices. Then the tensor product of the congruence
lattices, $\Con A$ and $\Con B$, is isomorphic to the congruence lattice of the
tensor product $A \otimes B$, in formula,
 \[
   \Con A \otimes \Con B \iso \Con (A \otimes B).
 \]
 \end{all}

 There are some stronger results stated in \cite{GLQu81}, see
Section~\ref{S:Discussion} for a discussion.

 In this paper, we generalize the main result of \cite{GLQu81} to infinite
lattices.  First, one has to observe that the isomorphism of the Main Result of
\cite{GLQu81} cannot be expected to hold for infinite lattices; indeed, easy
examples show (see Example~\ref{E:NonCon}) that $\Con (A \otimes B)$ is, in
general, very large when compared to $\Con A \otimes \Con B$. Indeed, the proof
of the main result of \cite{GLQu81} computes principal congruences of $A
\otimes
B$ in terms of the principal congruences of $A$ and $B$.  So the proper
generalization ought to change the congruence lattice to the semilattice with
zero of compact congruences. We denote by $\Conc L$ the semilattice with
zero of
compact congruences of the lattice~$L$.

 To state our result, we need the concept of capping.

 Let us call a subset $I$ of $A \times B$ a \emph{bi-ideal}, if $I$
contains $(A \times \set{0})  \uu  (\set{0} \times B)$, it is hereditary, and
it is join-closed in the sense that if $\vv<a_0, b>$, $\vv<a_1, b>\in I$, then
$\vv<a_0 \jj a_1, b> \in I$, and symmetrically. Then $A \otimes B$ can be
represented as the join-semilattice with zero of finitely generated bi-ideals
of $A \times B$.

 A bi-ideal $I$ is \emph{capped}, if there is a finite subset $C$ of $A \times
B$ such that $I$ is the hereditary subset of $A \times B$ generated by $C$
along
with $(A \times \set{0}) \uu (\set{0} \times B)$.  A tensor product $A \otimes
B$ is \emph{capped}, if all bi-ideals in the representation of $A \otimes
B$ are
capped.  If $A \otimes B$ is a capped tensor product, then it is a lattice.

 \begin{all}{Main Theorem}
 Let $A$ and~$B$ be lattices with zero.
 \begin{enumerate}
 \item If $A \otimes B$ is a lattice, then there is a natural embedding of
$\Conc A
\otimes \Conc B$ into $\Conc (A \otimes B)$.
 \item If $A \otimes B$ is a capped tensor product, then
 \[
   \Conc A \otimes \Conc B \iso \Conc (A \otimes B).
 \]
 \end{enumerate}
  \end{all}

 Most results in this paper are stated for the more general constructions
called
\emph{sub-tensor product} and \emph{capped sub-tensor product}, introduced in
Section~\ref{S:Sub-tensor}. There is a good reason for this, although, it
is not
evident in this paper.  In \cite{GrWe2}, we introduce a variant of the tensor
product construction, which we name \emph{box product}. The most important
advantage of the new construction is that the box product of two lattices is
always a lattice.  In \cite{GrWe2}, we state the analogue of the Main Theorem
for box products; it turns out that the proof for box products is very similar
to the proof in this paper.  So we introduce here sub-tensor products and
capped sub-tensor products, which serve as a common platform to prove the
results that apply in this paper and also in \cite{GrWe2}.

 In Section~\ref{S:Tensor}, we introduce the basic concepts and restate the
basic results on tensor products we shall need in this paper.

 The new concepts of \emph{L-homomorphism} and \emph{L-congruence} are
introduced in Section~\ref{S:L-congruences}.  We establish that L-homomorphisms
and L-congruences of join-semilattices with zero behave very much like
homomorphisms and congruences of lattices.  These concepts allow us to develop
results for semilattices that otherwise could only be obtained for lattices.
The main result in this section is Lemma~\ref{L:CanHom}, which lifts
L-homomorphisms to tensor products, generalizing a result---Lemma 3.17 of
\cite{GrWe2}---from finite lattices to arbitrary semilattices with zero.

 Sub-tensor products are introduced in Section~\ref{S:Sub-tensor}, where we
prove some basic properties.  In Section~\ref{S:TensProd}, if $A$ and $B$ are
lattices with zero and $C$ is a sub-tensor product of $A$ and $B$, then for a
compact congruence $\ga$ of $A$ and a compact congruence $\gb$ of $B$, we
define a congruence $\ga \odot_C \gb$ of $C$ and prove that there is a unique
\jz-homomorphism $\ge_C$ from $\Conc A \otimes \Conc B$ to $\Conc C$ such that,
for all $\ga \in \Conc A$ and all $\gb \in \Conc B$, we have $\ge_C(\ga \otimes
\gb) = \ga \odc{C} \gb$. In Section~\ref{S:Embedding}, we prove the Embedding
Theorem, claiming that this map $\ge_C$ is, in fact, an embedding; this
verifies the sub-tensor product version of the first statement of the Main
Theorem.

 In Section~\ref{S:capped}, we introduce capped sub-tensor products, and we
prove the Isomorphism Theorem that corresponds to the second statement of the
Main Theorem.

 In Section~\ref{S:Discussion}, we apply the Embedding Theorem and the
Isomorphism Theorem to get the two statements of the Main Theorem.  We also
discuss related results, in particular, some results from J. D. Farley
\cite{Farl96}, G. Gr\"atzer, H. Lakser, and R. W. Quackenbush \cite{GLQu81},
and G. Gr\"atzer and E. T. Schmidt \cite{GrSc94}, and state a number of open
problems.

 \section{Tensor products}\label{S:Tensor}
 Let $\E{S}_0$ denote the category of join-semi\-lat\-tices with zero with
join- and zero-preserving homomorphisms (\jz-semilattices with
\jz-homomor\-phisms).  For
$S \in \E S_0$, let $\go_S$ and $\gi_S$ denote the smallest and largest
congruence
of the \jz-semilattice $S$, respectively.

 Let $(A_i \mid i\in I)$ be a family of join-semilattices with $0$; the
\emph{direct sum} of this family in $\E S_0$, denoted by $\bigoplus(A_i
\mid i\in
I)$, is the \jz-semilattice of all $\vv<a_i \mid i\in I>$ such that $a_i
\in A_i$,
for $i \in I$, and $\setm{i\in I}{a_i\ne 0}$ is finite (the zero is the ``zero
vector'' and the join is formed componentwise). In fact, $\bigoplus(A_i
\mid i\in
I)$ is the coproduct of the family $(A_i\mid i\in I)$ in $\E S_0$.

 Now we introduce the basic definitions for this paper. In \cite{GLQu81}, the
tensor product of the objects $A$, $B \in \E{S}_0$ consists of certain
hereditary
subsets $X$ of $(A - \{0\}) \times(B - \{0\})$. In this paper, we find it more
convenient to consider hereditary subsets $X$ of $A \times B$ that contain
the set
 \[
   \nabla_{A,B} = (A \times \set{0})  \uu  (\set{0} \times B).
 \]

 We shall use a partial binary operation on $A \times B$: let $\vv<a_0,b_0>$,
$\vv<a_1,b_1>\in A \times B$; the \emph{lateral join} of $\vv<a_0,b_0>$ and
$\vv<a_1,b_1>$ is defined if $a_0 = a_1$ or $b_0 = b_1$, in which case, it
is the join, $\vv<a_0\jj a_1,b_0\jj b_1>$.

 \begin{definition}\label{D:etp}
 Let $A$ and $B$ be \jz-semilattices. A nonempty subset $I$ of $A \times B$
is a
\emph{bi-ideal} of $A \times  B$, if it satisfies the following conditions:
 \begin{enumerate}
 \item $I$ is hereditary;
 \item $I$ contains $\nabla_{A,B}$;
 \item $I$ is closed under lateral joins.
 \end{enumerate}
 The \emph{extended tensor product} of $A$ and $B$, denoted by $A \ootimes
B$, is
the lattice of all bi-ideals of $A \times B$.
 \end{definition}

 It is easy to see that $A \ootimes  B$, is an algebraic lattice. For $a
\in A$ and
$b \in B$, we define $a \otimes b \in A
 \ootimes  B$ by
 \[
   a \otimes b =\nabla_{A,B}  \uu  \setm{\vv<x, y> \in A \times B}{\vv<x,
y>  \leq
\vv<a, b>}
 \]
 and call $a \otimes b$ a \emph{pure tensor}.  A pure tensor is a
(one-generated)
principal bi-ideal. Let us call $(a_0 \otimes b_0) \jj (a_1 \otimes b_1)$ a
\emph{mixed tensor}, if $a_0  \leq  a_1$ and $b_0 \geq b_1$ or $a_0 \geq
a_1$ and
$b_0  \leq  b_1$.  A mixed tensor is a special type of join of two pure
tensors, a two-generated bi-ideal.

 \begin{definition}\label{D:tensor_product}
 Let $A$ and $B$ be \jz-semilattices. The \emph{tensor product} $A \otimes
B$ is
the \jz-sub\-semi\-lat\-tice of compact elements of $A  \ootimes  B$;
equivalently,
$A \otimes B$ is the \jz-sub\-semi\-lat\-tice generated in $A \ootimes  B$
by the
pure tensors.
 \end{definition}

 Since pure tensors and mixed tensors are compact bi-ideals, we conclude the
following:

 \begin{proposition}\label{P:puremixed}
 All pure tensors and all mixed tensors are elements of the tensor product.
 \end{proposition}

 The tensor product of \jz-semilattices may be a lattice.

 \begin{proposition}\label{P:lattice}
 $A \otimes B$ is a lattice if and only if it is closed under finite
intersection.
 \end{proposition}

 \begin{proof}
 Indeed, let $A \otimes B$ be a lattice and let $I \mm  J  =  K$, where
$I$, $J$,
and $K$ are compact bi-ideals.  If $K$ is not $I \ii J$, then $K \sci I \ii
J$ and
so there is a compact bi-ideal $H$ satisfying $K \sci H \sci I \ii J$,
contradicting that $I \mm J  =  K$.  The converse is trivial.
 \end{proof}

 This proposition is really a statement that holds for any algebraic lattice,
viewed as the ideal lattice of a \jz-semilattice.

 \begin{corollary}\label{C:twofinite}
 The tensor product of the finite lattices $A$ and $B$ is always a lattice.
 \end{corollary}

 \begin{proof}
 Indeed, then, $A \otimes B = A \ootimes B$ and $A \ootimes B$ is closed under
finite intersection, therefore, so is $A \otimes B$, and thus $A \otimes B$
is a
lattice.
 \end{proof}

 We shall next characterize the tensor product as a universal construction with
respect to bimorphisms.  So first we give the definition of a bimorphism.

 \begin{definition}\label{D:bimorphism}
 Let $A$, $B$, and $C$ be \jz-semilattices. A \emph{bimorphism} from $A
\times B$
to $C$ is a map $f \colon A \times B \to C$ such that
 \begin{enumerate}
 \item for all $\vv<a,b>\in A\times B$, $f(\vv<a,0>) = f(\vv<0,b>) = 0$;
\item for
all $a_0$, $a_1\in A$ and all $b\in B$,
 \[
   f(\vv<a_0 \jj a_1, b>) = f(\vv<a_0, b>) \jj f(\vv<a_1, b>);
 \] \item for all $a\in A$ and all $b_0$, $b_1\in B$,
 \[
   f(\vv<a, b_0 \jj b_1>) = f(\vv<a, b_0>) \jj f(\vv<a, b_1>).
 \]
 \end{enumerate}
 \end{definition}

 \begin{corollary}
 A bimorphism is an isotone map.
 \end{corollary}

 \begin{proof}
 Indeed, if $f \colon A \times B \to C$ is a bimorphism and $\vv<a_0,b_0>$,
$\vv<a_1,b_1> \in A \times B$ satisfy $\vv<a_0,b_0> \leq \vv<a_1,b_1>$, then
$f(\vv<a_0,b_0>) \jj f(\vv<a_0,b_1>) = f(\vv<a_0,b_1>)$, by
Definition~\ref{D:bimorphism}(iii), and so $f(\vv<a_0,b_0>) \leq
f(\vv<a_0,b_1>)$;
similarly, $f(\vv<a_0,b_1>)  \leq f(\vv<a_1,b_1>)$, by
Definition~\ref{D:bimorphism}(ii), from which the statement follows.
 \end{proof}

 \begin{proposition}\label{P:TensUniv}
 Let $A$ and $B$ be \jz-semilattices. Consider the map $\otimes \colon A
\times B
\to A \otimes B$ defined by $\vv<a, b> \mapsto a \otimes b$. Then $\otimes$
is a
universal bimorphism, that is, for every $C \in \E{S}_0$ and every
bimorphism $f
\colon A \times B \to C$, there exists a unique \jz-homomorphism $g \colon
A \otimes
B \to C$ such that $g(a\otimes b) = f(\vv<a, b>)$, for all $a \in A$ and $b
\in B$.
 \end{proposition}

 \begin{note}
 We could have defined $A \otimes B$ as the \jz-semilattice freely
generated by all
elements of $A \times B$, subjected to the relations
$\vv<a,0>=\vv<0,b>=0$,
$\vv<a_0\jj a_1, b> = \vv<a_0,b> \jj \vv<a_1,b>$ and
$\vv<a,b_0 \jj b_1> = \vv<a,b_0> \jj \vv<a,b_1>$,
for $a$, $a_0$, $a_1 \in A$ and $b$, $b_0$, $b_1 \in B$. With this definition,
Proposition~\ref{P:TensUniv} is evident since $A \otimes B$ is a free object.
However, for most computations, we need the representation of the elements
of $A
\otimes B$ by compact bi-ideals.  So if we define $A \otimes B$ as a free
object,
then we would replace Proposition~\ref{P:TensUniv} by the representation of $A
\otimes B$ as the compact bi-ideal \jz-semilattice of $A \times B$.
 \end{note}

 \begin{proof}
 It is routine to verify that the map $\otimes$ is a bimorphism. Since the pure
tensors generate $A\otimes B$ as a \jz-semilattice, the uniqueness statement is
also trivial.

 For a given bimorphism
$f \colon A \times B \to C$, we now prove the existence of $g$ such
that $g(a \otimes b) = f(\vv<a, b>)$, for all $a \in A$ and $b \in B$.  Let
$\E D$
be the set of all subsets $X$ of $A \times B$ such that
$\JJm{f(\vv<x,y>)}{\vv<x,y>\in X}$ is defined in $C$. For every $X\in\E D$, put
 \[
    h(X) = \JJm{f(\vv<x,y>)}{\vv<x,y>\in X}
 \]
 (the join is formed in $C$).

\setcounter{claim}{0}
 \begin{claim}\label{C:1}
 For every $X \in \E D$, the bi-ideal $\ol X$ of $A \times B$ generated by
$X$ also
belongs to $\E D$ and $h(X)  =  h(\ol X)$.
 \end{claim}

 \begin{proof}
 Let $X_0$ be $X  \uu  \nabla_{A, B}$, and, for every integer $n > 0$, let
$X_{n}$
be the hereditary set generated by lateral joins of elements of $X_{n-1}$.
Obviously, $X_n \in \E D$ with $h(X) = h(X_n)$, for all $n \geq 0$.  Since
 \[
   \ol X =  \UUm{X_n}{n \geq 0},
 \]
 the statement follows.
 \end{proof}

 The proof of the following claim is obvious:

\begin{claim}\label{C:2}
 The set $\E D$ is closed under finite unions, and $h$ is a
\jz-homomorphism from
$\vv<\E{D}; \uu,\nabla_{A,B}>$ to $\vv<C;\jj,0>$.
 \end{claim}

 Since $h$ is defined on all pure tensors, it follows from Claims \ref{C:1} and
\ref{C:2} that $h$ is defined on $A \otimes B$, and that the restriction
$g$ of $h$
to $A \otimes B$ is a \jz-homomorphism from $A \otimes B$ to $C$. For all
$\vv<a,
b> \in A \times B$, it is obvious that $g(a \otimes b) = f(\vv<a, b>)$.
 \end{proof}

 This characterization of the universal bimorphism on $A \times B$ shows that
$\otimes$ defines, in fact, a \emph{bifunctor} on $\E S_0$. This allows us to
prove, in a routine fashion, the two following basic categorical results.

\begin{proposition}\label{P:AsCom}
 The tensor product operation is associative and commutative. Thus, if $A$,
$B$,
and $C$ are \jz-semilattices, then the following isomorphisms hold:
 \begin{align*}\label{Eq:AsCom}
   (A \otimes B)\otimes C &\iso A \otimes(B \otimes C);\\
              A \otimes B &\iso B \otimes A.
 \end{align*}
 \end{proposition}
 Note that these isomorphisms are \emph{natural} in the categorical sense.
We leave
the details to the reader.

\begin{proposition}\label{P:BasPres}
 Let $B$ be a \jz-semilattice. Then the functor
 \[
   \ul{\phantom{X}}\otimes B
 \]
 preserves direct sums and directed colimits in $\E S_0$. In particular, if
$(A_i
\mid i\in I)$ is a family of \jz-semilattices, then
 \begin{equation*}
 \bigoplus(A_i\mid i \in I) \otimes B \iso
 \bigoplus(A_i \otimes B \mid i\in I).
 \end{equation*}

 Similarly, if $I$ is a directed set and $A = \varinjlim_i A_i$ with
respect to a
limit system on $(A_i\mid i\in I)$, then
 \begin{equation*}
   \varinjlim (A_i \mid i)\otimes B \iso \varinjlim(A_i\otimes B
     \mid i).
 \end{equation*}
 \end{proposition}

 The following purely arithmetical formulas are due to G.~A. Fraser
\cite{Fras76}.

\begin{lemma}\label{L:IntersTens}
 Let $A$ and $B$ be \jz-semilattices.  Let $a_0$, $a_1 \in A$ and $b_0$,
$b_1 \in
B$ such that $a_0 \mm  a_1$ and $b_0 \mm  b_1$ both exist.
\begin{enumerate}
\item The intersection of two pure tensors is a pure tensor, in fact,
 \begin{equation*}
   (a_0 \otimes b_0) \ii (a_1 \otimes b_1)
     = (a_0 \mm  a_1) \otimes (b_0 \mm  b_1).
 \end{equation*}
 \item The join of two pure tensors is the union of four pure tensors, in fact,
 \begin{align*}
      (a_0 \otimes b_0) \jj (a_1 \otimes b_1) & = \\
    (a_0 \otimes b_0)  \uu  (a_1 \otimes b_1)
     & \uu  ((a_0 \jj a_1) \otimes (b_0 \mm  b_1))
  \uu  ((a_0 \mm  a_1) \otimes (b_0 \jj b_1)).
 \end{align*}
 \item A mixed tensor is a union of two pure tensors, that is, if $a_0
\leq  a_1$
and $b_0 \geq b_1$, or $a_0 \geq a_1$ and $b_0 \leq  b_1$, then
\begin{equation*}
   (a_0 \otimes b_0) \jj (a_1 \otimes b_1) =
   (a_0 \otimes b_0) \uu (a_1\otimes b_1).
 \end{equation*}
 \item Let $A$ and $B$ be lattices with zero.  Then
 \begin{gather*}
  \JJm{a_i \otimes b_i}{i < n} \mm \JJm{c_j \otimes d_j}{j < m}\\
     = \UU (p(a_1, \dots, a_{n - 1}) \mm q(c_1, \dots, c_{m - 1}))
      \otimes (p^{\mrm{d}}(b_1, \dots, b_{n - 1}) \mm
      q^{\mrm{d}}(d_1, \dots, d_{m - 1})),
 \end{gather*}
 where $p^{\mrm{d}}$ and $q^{\mrm{d}}$ are the duals of $p$ and $q$,
respectively,
and where the union is for all $p \in F(n)$ and $q \in F(m)$.
 \end{enumerate}
 \end{lemma}

 \begin{corollary}\label{C:distr}
 Let $A$ and $B$ be lattices with zero.  Let $a_0$, $a_1 \in A$ and $b_0$,
$b_1 \in
B$ satisfy $a_0  \leq  a_1$ and $b_0 \geq b_1$, or $a_0 \geq a_1$ and $b_0
\leq
b_1$.  Set $I = a_0 \otimes b_0$ and $J = a_1 \otimes b_1$.  Then the
distributive
law
 \[
   (I \jj J) \mm H = (I \mm H) \jj (J \mm H)
 \]
 holds, for any $H \in A \otimes B$.
 \end{corollary}

We can rephrase the statements of this lemma with the following concept:

 \begin{definition}\label{D:capping}
 Let $I$ be a bi-ideal of $A \times B$.  A \emph{capping} of $I$ is a
\emph{finite}
subset $C$ of $A \times B$ so that
 \[
   I = \setm{x \in A \times B}{x \leq i, \text{ for some } i \in C} \uu
\nabla_{A,B},
 \]
 that is, $I$ is the hereditary set generated by $C$ in $A \times B$ along with
$\nabla_{A,B}$.  A \emph{capped bi-ideal} is a bi-ideal with capping.
 \end{definition}

 For instance, $a \otimes b$ is capped by $\set{\vv<a, b>}$ and
Lemma~\ref{L:IntersTens}(iii) can be restated as follows: a mixed tensor $(a_0
\otimes b_0) \jj (a_1 \otimes b_1)$ (where $a_0 \leq a_1$ and $b_0 \geq
b_1$, or
$a_0 \geq a_1$ and $b_0  \leq b_1$) is capped by $\set{\vv<a_0, b_0>, \vv<a_1,
b_1>}$.

 A capped bi-ideal is compact, but, in general, a compact bi-ideal may not be
capped. In \cite{GrWe1}, the reader may find examples of compact bi-ideals
that are not capped.  For instance, let $a$, $b$, and $c$ be the atoms of
$M_3$, let $x_0$, $x_1$, $x_2$ be the free generators of $\tup{F}(3)$, the
free lattice on three generators; then $M_3 \otimes \tup{F}(3)$
contains such examples, for instance, the bi-ideal generated by
$\set{\vv<a, x_0>, \vv<b, x_1>, \vv<c, x_2>}$. See Section~\ref{S:capped}, for
applications of this concept to tensor products.

 \section{L-congruences}\label{S:L-congruences}

 In this section, we introduce L-homomorphisms and L-congruences. These
concepts
allow us to develop results for \jz-semilattices that otherwise could only be
obtained for lattices.

\begin{definition}
 Let $A$ and $B$ be \jz-semilattices. Let $f \colon A \to B$ be a
\jz-homomorphism.
We shall say that $f$ is an \emph{L-homomorphism}, if for all $a_0$, $a_1
\in A$
and $b \in B$,
 \[
   b  \leq  f(a_0)\text{\q and\q }b  \leq  f(a_1)
 \] imply the existence of an $x \in A$ such that
 \[
   x  \leq  a_0,\q x  \leq  a_1\text{,\q and\q }b  \leq  f(x).
 \]
 An \emph{L-congruence} of a \jz-semilattice $A$ is the kernel of an
L-homo\-morph\-ism from $A$ to some \jz-semilattice $B$. \end{definition}

In this definition, the \emph{kernel} of a map is the equivalence relation
induced
by it.  The kernel of $f$ will be denoted by $\ker f$.

 \begin{corollary}\label{C:LcongSpecial}
 If $f \colon A \to B$ is an L-homomorphism, $a_0$, $a_1 \in A$, and
$f(a_0) \leq
f(a_1)$, then there is an $\ol a_0 \in A$ with $\ol a_0 \leq a_0$ and $\ol
a_0 \leq
a_1$ such that $f(a_0) = f(\ol a_0)$.
 \end{corollary}

\begin{proof}
 Choose $b = f(a_0)$.  Then $b  \leq  f(a_0)$ and $b \leq f(a_1)$, so there
is an
$\ol a_0 \in B$ satisfying  $\ol a_0 \leq a_0$, $\ol a_0 \leq  a_1$, and $b
\leq
f(\ol a_0)$. Since $f(\ol a_0) \leq f(a_0) = b$, we conclude that $f(\ol a_0) =
f(a_0)$.
 \end{proof}

 \begin{proposition}\label{P:new}
 Let $A$ and $B$ be \jz-semilattices, and let $f\colon A\to B$ be a
L-homomorphism. Then the following holds:
\begin{enumerate}
\item $f$ is a \emph{partial meet-homomorphism}, that is, for all $n>0$
and all $a_0$, \dots, $a_{n-1}\in A$, if $a=a_0\mm\cdots\mm a_{n-1}$ exists in
$A$, then $b=f(a_0)\mm\cdots\mm f(a_{n-1})$ exists in $B$, and $b=f(a)$.

\item If, in addition, $f$ is one-to-one, then $f$ is a \emph{partial
meet-embedding}, that is, for all $n>0$ and all $a_0$, \dots, $a_{n-1}\in A$,
$a=a_0\mm\cdots\mm a_{n-1}$ exists in $A$ iff $b=f(a_0)\mm\cdots\mm f(a_{n-1})$
exists in $B$, and then, $b=f(a)$.
\end{enumerate}

Conversely, if $A$ and $B$ are lattices, then any
lattice homomorphism from $A$ to $B$ is an L-homomorphism and any
L-congruence of
$A$ is a lattice congruence of $A$.
 \end{proposition}

 \begin{proof}
(i) If $A$ and $B$ are \jz-semilattices and $f \colon A \to B$ is an
L-homomorphism, we prove that $f$ is a partial meet-homomorphism. Let $n>0$,
let $a_0$, \dots, $a_{n - 1}\in A$, and let
$a=a_0\mm\cdots\mm a_{n-1}$ be defined in $A$. Since $f$ is isotone,
$f(a)\leq f(a_0)$, \dots, $f(a_{n-1})$ in $B$. Conversely, let $b \leq f(a_0)$,
\dots, $f(a_{n-1})$ in $B$. Since $f$ is an L-homomorphism, there exists
$x\in A$ such that $x\leq a_i$, for all $i<n$, and $b\leq f(x)$. Since
$x\leq a$ and $f$ is isotone, we
have $b \leq f(x) \leq f(a)$. This proves that
$f(a_0)\mm\cdots\mm f(a_{n-1})$ is defined and equals $f(a)$.

(ii) Now, suppose that $f$ is one-to-one. Thus, since $f$ is a
join-homomorphism,
$f$ is an order-embedding, that is, $f(x)\leq f(y)$ iff $x\leq y$, for all $x$,
$y\in A$. Now let $n>0$, let $a_0$, \dots, $a_{n-1}\in A$. Suppose that
$b=f(a_0)\mm\cdots\mm f(a_{n-1})$ exists in $B$. Since $b\leq f(a_i)$ for all
$i$ and since $f$ is an L-homomorphism, there exists $a\in A$ such that
$a\leq a_i$ for all $i$, and $b\leq f(a)$. Since $f$ is isotone,
$f(a)\leq f(a_i)$ for all $i$, thus $f(a)\leq b$; so $b=f(a)$.
For all $x\in A$ such that $x\leq a_i$ for all $i$, we have $f(x)\leq
f(a_i)$ for
all $i$, thus $f(x)\leq b$, that is, $f(x)\leq f(a)$; whence $x\leq a$. This
proves that $a=a_0\mm\cdots\mm a_{n-1}$. Thus (ii) holds as well.

Now let $A$ and $B$ be lattices with zero and
let $f\colon A \to B$ be a lattice homomorphism. We can choose
$x = a_0 \mm a_1$
(in the definition of L-homomorphism) to verify that $f$ is an L-homomorphism.

Finally, the result about L-homomorphisms implies immediately the result about
L-congruences.
 \end{proof}

Another connection between L-homomorphisms and lattice homomorphisms is the
following:

 \begin{proposition}\label{P:idealhom}
 Let $A$ and $B$ be \jz-semilattices and let $f\colon A\to B$ be an
L-homomorphism. For an ideal $I$ of $A$, define $\ol{f}(I)$ as the
hereditary subset of $B$ generated by $f(I)$. Then $\ol{f}$ is a join-complete,
$\{\jj,\mm,0\}$ homomorphism of $\Id A$ to $\Id B$, and it has the property
that
the image of a principal ideal is principal.

Conversely, let $g\colon\Id A\to\Id B$ be a join-complete, $\{\jj,\mm,0\}$
homomorphism with the property that the image of a principal ideal is
principal.
Then there exists a unique L-homomorphism $f\colon A\to B$ such that
$g=\ol{f}$.
 \end{proposition}

 We leave the proof to the reader. All the steps in this proof are easy; we use
that $f$ is an L-homomorphism in verifying that $\ol f(I) \ii \ol f(J) \ci
\ol f(I
\ii J)$.

 For a set $X$  and a binary relation $\ga$ on $X$, we use the notation $x
\equiv_\ga y$, for $\vv<x, y>\in\ga$.

 Let $A$ and $B$ be \jz-semilattices. Let $\ga$ be an L-congruence of $A$
and let
$\gb$ be an L-congruence of $B$. Then $\ga \times \gb$ is a congruence on
$A \times
B$.  Now we define how $\ga \times \gb$ can be naturally extended to a
congruence
$\ga \congtimes \gb$ of $A\otimes B$.

 \begin{definition}
 Let $A$ and $B$ be \jz-semilattices. Let $\ga$ be an L-congruence of $A$
and let
$\gb$ be an L-congruence of $B$.  Define a binary relation $\ga \bboxtimes
\gb$ on
$A  \ootimes  B$ as follows: for $H$, $K \in A \ootimes  B$, let
$H\equiv_{\ga\bboxtimes\gb}K$ iff, for all $\vv<x, y>\in H$, there exists an
$\vv<x', y'>\in K$ such that $x \equiv_{\ga} x'$ and $y \equiv_{\gb} y'$, and
symmetrically.

 Let $\ga \congtimes \gb$ be the restriction of $\ga \bboxtimes\gb$ to
$A\otimes B$.
 \end{definition}

 The following result allows us to give a useful explicit description of
the effect
of the tensor product bifunctor $\otimes$ on two homomorphisms in~$\E S_0$.

\begin{lemma}\label{L:CanHom}
 Let $A$, $A'$, $B$, $B'$ be \jz-semilattices, let $f \colon A \to A'$ and $g
\colon B \to B'$ be L-homomorphisms. For a bi-ideal $I$ of $A \times B$, define
$h(I)$ as the the hereditary subset of $A' \times B'$ generated by the
image of $I$
under $f \times g$ with $\nabla_{A', B'}$, that is,
 \[
   h(I)  =  \nabla_{A', B'}  \uu  \setm{\vv<u, v> \in A' \times B'}
 {u  \leq  f(x) \text{ and } v \leq g(y) \text{, for some } \vv<x, y> \in I}.
 \]
 Then the following properties hold:
 \begin{enumerate}
 \item $h(I)$ is a bi-ideal of $A'\times B'$.
 \item The map $h$ is a lattice homomorphism from $A \ootimes B$ to $A'
\ootimes
B'$.
 \item $h(A \otimes B) \ci A' \otimes B'$ and the restriction of $h$ from $A
\otimes B$ to $A'\otimes B'$ equals $f\otimes g$.
 \item $f \otimes g$ is an L-homomorphism.
 \item $\ker h = \ker f \bboxtimes \ker g$. Thus, $\ker(f \otimes g) = \ker f
\congtimes \ker g$.
 \end{enumerate}
 \end{lemma}

 \begin{proof}\hfill

 (i) By definition, $h(I)$ is hereditary and contains $\nabla_{A', B'}$. Let
$\vv<x'_0, y'>$, $\vv<x'_1, y'>\in h(I)$; we prove that $\vv<x'_0 \jj x'_1,
y'> \in
h(I)$.  If $\vv<x'_0, y'> \in A' \times \set{0}$, then $y' = 0$ and
$\vv<x'_0 \jj
x'_1, 0> \in h(I)$ by the definition of $h(I)$. Similarly, the conclusion is
obvious if $\vv<x'_0, y'> \in \set{0} \times B'$.  So assume that
$\vv<x'_0, y'>$,
$\vv<x'_1, y'>\nin \nabla_{A',B'}$. Then by the definition of $h(I)$, there
exist
$\vv<x_0, y_0>$, $\vv<x_1, y_1>\in I$ such that $\vv<x'_0, y'> \leq
\vv<f(x_0),
g(y_0)>$ and $\vv<x'_1, y'> \leq  \vv<f(x_1), g(y_1)>$. Since $g$ is an
L-homomorphism, there exists a $y\in B$ such that
 \[
   y  \leq  y_0,\ y  \leq  y_1 \text{, and } y'  \leq  g(y).
 \]
 Since $I$ is hereditary, it follows that $\vv<x_0, y>$, $\vv<x_1, y> \in
I$, and
so $\vv<x, y>\in I$ with $x  =  x_0 \jj x_1$. Since $\vv<x'_0 \jj x'_1, y'>
\leq
\vv<f(x), g(y)>$, this proves that $\vv<x'_0 \jj x'_1, y'> \in h(I)$. By
symmetry,
this proves that $h(I)$ is a bi-ideal of $ A'\times B'$  (and so $h(I) \in A'
\ootimes  B'$).

 (ii) $h$ is a join-homomorphism.  Indeed, let $I$ and $J$ be bi-ideals of $A
\times B$. It is obvious that $h(I) \jj h(J) \ci h(I \jj J)$. Conversely, the
following set
 \[
    X = \setm{\vv<x, y> \in A \times B}{\vv<f(x), g(y)> \in h(I)
      \jj h(J)}
 \]
 is a bi-ideal, and it obviously contains $I$ and $J$. Thus it contains $I
\jj J$,
which implies that $h(I \jj J) \ci h(I)\jj h(J)$.

 Now we prove that $h(I \ii J)  =  h(I) \ii h(J)$, for bi-ideals $I$ and
$J$ of $A
\times B$. To prove the nontrivial containment, let $\vv<u, v> \in h(I) \ii
h(J)$,
and we want to prove that $\vv<u, v> \in h(I \ii J)$.  This is trivial if $u =
0_{A'}$ or $v = 0_{B'}$.  So assume that $u$ and $v$ are nonzero.  Then, by
definition, there are $\vv<x', y'> \in I$ and $\vv<x'', y''>\in J$ such that
$\vv<u, v>  \leq \vv<f(x'), g(y')>$ and $\vv<u, v> \leq  \vv<f(x''), g(y'')>$.
Since $f$ and $g$ are L-homomorphisms, there are $x \in A$ and $y \in B$
such that
$x  \leq  x'$, $x''$ and $y  \leq  y'$, $y''$ and $\vv<u, v>  \leq  \vv<f(x),
g(y)>$. Since $\vv<x, y> \in I \ii J$, we have proved that $\vv<u, v> \in
h(I \ii
J)$.

 (iii) is obvious, because $h$ is a join-homomorphism and, for all $\vv<a,
b> \in A
\times B$, we have that $h(a \otimes b) = f(a) \otimes g(b)$ (we use here
the fact
that both $f$ and $g$ are zero-preserving).

 (iv) Put $h'=f \otimes g$. Let $I_0$, $I_1 \in A \otimes B$ and $J \in A'
\otimes
B'$ such that $J \leq h'(I_0)$, $h'(I_1)$. Since $h$ is a lattice
homomorphism, we
obtain that $J \leq  h(I)$ with $I = I_0 \ii I_1$. Obviously, $h$ is a complete
join-homomorphism and $I$ is the directed union of all compact bi-ideals of $A
\times B$ contained in $I$, therefore, there exists a compact bi-ideal $I'$
of $A
\times B$ such that $I' \leq I$ and $J \leq h'(I')$. This proves that $h'$
is an
L-homomorphism.

 (v) It suffices to prove the first statement. To prove that $\ker h \leq
\ker f
\bboxtimes \ker g$, take $I$, $J\in A \ootimes B$ and assume that $I
\equiv_{\ker
h} J$, that is, $h(I) = h(J)$.  We prove that for every $\vv<x, y> \in I$
there is
$\vv<\ol x, \ol y> \in J$ such that $x \equiv_{\ker f} \ol x$ and $y
\equiv_{\ker g}
\ol y$.  If $f(x) = 0_{A'}$, then $x \equiv_{\ker f} 0_A$ and $\vv<0_A,
y>\in J$,
so $\ol x = 0_A$ and $\ol y = y$ will do. Argue similarly for $g(y) =
0_{B'}$. If
both $f(x)$ and $g(y)$ are nonzero, then $\vv<f(x), g(y)>$ is majorized by some
$\vv<f(\ol x), g(\ol y)>$,  where $\vv<\ol x, \ol y> \in J$. Since $f$ is an
L-homomorphism, just as in Corollary~\ref{C:LcongSpecial}, there is an $x_0
\in A$
such that $f(x) = f(x_0)$ and $x_0 \leq x$, $\ol x$.  Similarly, there is
an $y_0
\in B$ such that $f(y) = f(y_0)$ and $y_0 \leq y$, $\ol y$.  Since $\vv<\ol
x, \ol
y> \in J$, it follows that $\vv<x_0, y_0> \in J$ and, obviously, $x
\equiv_{\ker f}
x_0$ and $y \equiv_{\ker g} y_0$. By symmetry, this proves that $I
\equiv_{\ker f
\bboxtimes \ker g} J$. The converse is easy.
 \end{proof}

 The results of Lemma~\ref{L:CanHom} are formulated for L-homomorphisms. The
situation for general \jz-homomorphisms is quite different.  Let us call a
\jz-semilattice $S$ \emph{flat}, if for every \jz-semilattice embedding $f
\colon A
\hookrightarrow B$, the map
 \[
   \mathrm{id}_S \otimes f \colon S \otimes A \to S \otimes B
 \]
 is an embedding. We can prove that \emph{a \jz-semilattice is flat if and
only if
it is distributive}; see \cite{GrWe3}.

 \begin{corollary}\label{C:NabCon}
 Let $A$ and $B$ be \jz-semilattices, let $\ga$ be an L-congruence of $A$,
and let
$\gb$ be an L-congruence of $B$. Then the following properties hold:
 \begin{enumerate}
 \item The relation $\ga \congtimes \gb$ is an L-congruence of $A\otimes B$.
 \item If $A \otimes B$ is a lattice, then $\ga \congtimes \gb$ is a lattice
congruence on $A \otimes B$.
 \item Let $f$ (resp., $g$) be the canonical projection from $A$ onto $A/\ga$
(resp., from $B$ onto $B/\gb$). Then $f \otimes g$ factors into an
isomorphism from
$A \otimes B/\ga \congtimes \gb$ onto $(A/\ga) \otimes (B/\gb)$.
 \item If $A \otimes B$ is a lattice, then $(A/\ga) \otimes (B/\gb)$ is a
lattice.
 \end{enumerate}
 \end{corollary}
 \begin{proof}
 All these statements are obvious from the results of this section.
 \end{proof}

\begin{corollary}\label{C:Subl}
 Let $A$, $A'$, $B$, $B'$ be lattices with zero such that $A$ is a
$\set{0}$-sublattice of $A'$ and $B$ is a $\set{0}$-sublattice of $B'$. Let $f$
(resp., $g$) denote the inclusion map from $A$ into $A'$ (resp., $B$ into
$B'$).
Then $f\otimes g$ is a join-embedding of $A \otimes B$ into
$A' \otimes B'$, and it is a partial meet-embedding.

In particular, if $A' \otimes B'$ is a lattice, then $A \otimes B$ is a lattice
and $f \otimes g$ is a lattice embedding from $A \otimes B$ into $A'
\otimes B'$.
 \end{corollary}
 \begin{proof}
 By Lemma~\ref{L:CanHom}, $h=f\otimes g$ is a one-to-one L-homomorphism.
Thus, by
Proposition~\ref{P:new}, it is a partial meet-embedding.

If $A'\otimes B'$ is a lattice, we prove that $A\otimes B$ is a lattice.
Let $X$,
$Y\in A\otimes B$. Then $h(X)$, $h(Y)\in A'\otimes B'$, thus, since $A'\otimes B'$
is a lattice, $h(X)\mm h(Y)$ exists in $A'\otimes B'$. Since $h$ is a
partial
meet-embedding, $X\mm Y$ exists in $A\otimes B$; whence $A\otimes B$ is a
lattice. Since $h$ is a L-homomorphism from one lattice to the other, it is
also, by Proposition~\ref{P:new}, a lattice homomorphism.
 \end{proof}

 \begin{corollary}\label{C:Subl2}
 Let $A$, $A'$, $B$ and $B'$ be lattices with zero such that $A$ is a
sublattice of
$A'$ and $B$ is a sublattice of $B'$ (we assume neither $0_A = 0_{A'}$ nor
$0_B =
0_{B'}$). If $A' \otimes B'$ is a lattice, then $A \otimes B$ is a lattice.
 \end{corollary}

 \begin{proof}
 Put $A'' = A  \uu  \set{0_{A'}}$ and $B'' = B \uu  \set{0_{B'}}$. Then
$A''$ is a
$\set{0}$-sublattice of $A'$ and $B''$ is a $\set{0}$-sublattice of $B'$,
thus, by
Corollary~\ref{C:Subl}, $A'' \otimes B''$ is a lattice. Furthermore, $A$ is a
quotient of $A''$ (by the map that sends $0_{A'}$ to $0_A$ and all $x\in A$ to
$x$---a retraction). Similarly, $B$ is a quotient of $B''$. Therefore, by
Corollary~\ref{C:NabCon}(iv), $A \otimes B$ is a lattice.
  \end{proof}

 These results are related to Lemma 3.17 in \cite{GLQu81}, in which it is
proved
that if $A'$ is a finite lattice and the lattice $B'$ with $0$ is $A'$-lower
bounded (see Section~\ref{S:GLQu81}) and $A$, $B$ are $\set{0}$-sublattices
of $A'$
and $B'$, respectively, then $A \otimes B$ has a natural embedding into $A'
\otimes
B'$. Note that under these conditions, $A' \otimes B'$ is a lattice.

\section{Sub-tensor products}\label{S:Sub-tensor}

 \begin{definition}\label{D:TensLatt}
 Let $A$ and $B$ be lattices with zero. A \emph{sub-tensor product} of $A$
and $B$
is a subset $C$ of $A\otimes B$ satisfying the following conditions:
 \begin{enumerate}
 \item $C$ contains all the mixed tensors in $A\otimes B$;
 \item $C$ is closed under finite intersection;
 \item $C$ is a lattice with respect to containment.
 \end{enumerate}
 \end{definition}

 Note about this concept:
 \begin{enumerate}
 \item Every pure tensor $a \otimes b$ ($a \in A$, $b \in B$) belongs to
$C$ and
$0_{A \ootimes B} = \nabla_{A, B} \in C$.
 \item  $A \otimes B$ is not a meet-semi\-lattice (see \cite{GrWe1} for an
example).  That is why we require that $C$ be a meet-subsemilattice of $A
\ootimes
B$, not of $A \otimes B$.
 \item A sub-tensor product of $A$ and $B$ may not be a
\emph{join}-subsemilattice
of $A \otimes B$ (although it is a join-semilattice in its own right).
 \item Let $H_0$, \dots, $H_{n-1} \in C$. Then $\JJm{H_i}{i < n}$ in $A
\otimes B$
is, in general, smaller than $\JJm{H_i}{i < n}$ in $C$. Note, however,
Proposition~\ref{P:sub-tensor}(iv).  If we want to remind the reader that
the join is formed in $C$, we use the notation $\jj_C$ and $\JJ_C$.
 \end{enumerate}

 We now list some simple properties of sub-tensor products.

 \begin{proposition}\label{P:sub-tensor}
 Let $A$ and $B$ be lattices with zero and let $C$ be a sub-tensor product
of $A$
and $B$.  Let $a_0$, $a_1 \in A$ and $b_0$, $b_1 \in B$ and let $H$, $H_i
\in C$,
$i < n$.
 \begin{enumerate}
 \item If $H = \IIm{H_i}{i < n}$, then $H = \MMm{H_i}{i < n}$ in $C$.
 \item $(a_0 \otimes b_0) \mm (a_1 \otimes b_1)
       = (a_0 \mm  a_1) \otimes (b_0 \mm  b_1)$ in $C$.
 \item If $H = \JJm{H_i}{i < n}$ in $A \otimes B$, then
$H = \JJm{H_i}{i < n}$ in $C$.
 \item Every $H \in C$ can be represented in the form $H = \JJm{a_i \otimes
b_i}{i
< n}$ (the join formed in $C$), where $a_i \in A$ and $b_i \in B$, $i < n$.
 \item If $a_0  \leq  a_1$ and $b_0 \geq b_1$, or $a_0 \geq a_1$ and
$b_0 \leq b_1$.  Then
 \[
    (a_0 \otimes b_0) \jj (a_1 \otimes b_1)
       = (a_0 \otimes b_0) \uu (a_1 \otimes b_1)
 \]
 holds in $C$.
 \item Let $a_0  \leq  a_1$ and $b_0 \geq b_1$, or $a_0 \geq a_1$ and $b_0
\leq
b_1$.  Set $I = a_0 \otimes b_0$ and $J = a_1 \otimes b_1$.  Then $I$, $J
\in C$
and the distributive law
 \[
   (I \jj J) \mm H = (I \mm H) \jj (J \mm H)
 \]
 holds in $C$, for any $H \in C$.
 \end{enumerate}
 \end{proposition}

 \begin{proof}
 If $H = \IIm{H_i}{i < n}$, then $H = \MMm{H_i}{i < n}$ in the lattice $A
\ootimes
B$.  Since $C$ is a subposet of $A \ootimes B$, it follows that $H =
\MMm{H_i}{i <
n}$ in $C$, proving (i).

 By Lemma~\ref{L:IntersTens}(i), $(a_0 \otimes b_0) \ii (a_1 \otimes b_1)$
is $(a_0
\mm  a_1) \otimes (b_0 \mm  b_1)$ and this element is in $C$, by assumption, so
(ii) follows from (i).

 Let $H = \JJm{H_i}{i < n}$ in $A \otimes B$.  Since $C$ is a subposet of $A
\otimes B$, it follows that $H = \JJm{H_i}{i < n}$ in $C$, proving (iii).

 By the definition of $A \otimes B$ (Definition~\ref{D:tensor_product}),
every $H
\in C$ can be represented in the form $H = \JJm{a_i \otimes b_i}{i < n}$, where
$a_i \in A$ and $b_i \in B$, $i < n$ and the join is formed in $A \otimes
B$; so by
(iii), $H = \JJm{a_i \otimes b_i}{i < n}$ in $C$, proving (iv).

 (v) follows similarly from (iii).

 Finally, (vi) follows from (v).
 \end{proof}

 This section and the next two sections deal with the congruence structure of a
sub-tensor product $C$ of $A$ and $B$.  So it is reasonable to ask whether
there is
such a~$C$.  We show in \cite{GrWe2} that, for any lattices with zero $A$
and $B$,
there exists a sub-tensor product $C$ of $A$ and $B$.

 However, in this paper, the main result is in Section~\ref{S:capped}, where we
assume that $A \otimes B$ is capped (meaning that all the bi-ideals of $A
\times B$
are capped) and, therefore, a lattice.  In this case, we always have at
least one
sub-tensor product, namely, $A \otimes B$:

 \begin{proposition}\label{P:AotBTensLatt}
 $A \otimes B$ is a lattice if and only if it is a sub-tensor product of $A$
and~$B$.
 \end{proposition}

 \begin{proof}
 If $A \otimes B$ is a lattice, then \ref{D:TensLatt}(i) holds by
Proposition~\ref{P:puremixed}; \ref{D:TensLatt}(ii) holds by
Proposition~\ref{P:lattice}; while \ref{D:TensLatt}(iii) holds, by assumption.

  Conversely, if \ref{D:TensLatt}(i)--(iii) hold for $C = A \otimes B$, then $A
\otimes B$ is a lattice by Proposition~\ref{P:lattice}.
 \end{proof}

 \begin{proposition}\label{P:NiceArith}
 Let $A$ and $B$ be lattices with zero and let $C$ be a sub-tensor product
of $A$
and $B$.  The tensor product operation, viewed as a mapping from $A \times
B$ to
$C$, is a bimorphism.
 \end{proposition}

 \begin{proof}
 Let $a_0$, $a_1 \in A$ and $b \in B$. Since $(a_0 \jj a_1) \otimes b$
belongs to
$C$ and since it equals $(a_0 \otimes b) \jj (a_1 \otimes b)$ in $A \otimes
B$, by
Proposition~\ref{P:sub-tensor}(iii), the same holds in $C$. By symmetry, the
conclusion follows.
 \end{proof}

 For every lattice congruence $\ga$ of $A$ and $\gb$ of $B$, denote by
$\ga\tim{C}\gb$ the restriction of $\ga \congtimes \gb$ from $A \otimes B$
to $C$.
Then define
 \begin{align*}
   \ge_{A, C}(\ga) &= \ga\tim{C}\go_B,\\
   \ge_{B, C}(\gb) &= \go_A\tim{C}\gb. \end{align*}

 \begin{lemma}\label{L:e0e1hom}
 Let $A$ and $B$ be lattices with zero and let $C$ be a sub-tensor product
of $A$
and $B$. The map $\ge_{A, C} \colon \vv<\Con A; \jj, \go_A, \gi_A> \to \vv<\Con
C;\jj, \go_C, \gi_C>$ is a homomorphism (a $\set{\jj, 0, 1}$-homomorphism).
And,
similarly, for $\ge_{B, C}$.
 \end{lemma}

\begin{proof}
 It is clear that $\ge_{A, C}(\go_A) = \go_{C}$. Let $H$ be an element of
$C$ and
let $\vv<x,y>\in H$. Then $x \equiv_{\gi_A}0_A$ and $\vv<0_A, y>\in 0_{C}$
and so
$\vv<x, y> \equiv_{\gi_A \congtimes \go_B}\vv<0_A, y>$; this proves that $H
\equiv_{\ge_{A, C}(\gi_A)} 0_{C}$; whence $\ge_{A, C}(\gi_A) = \gi_{C}$.

 Now let $\ga_0$, $\ga_1 \in \Con A$. It is obvious that $\ge_{A, C}(\ga_0) \jj
\ge_{A, C}(\ga_1) \leq \ge_{A, C}(\ga_0 \jj \ga_1)$. Conversely, let $H$
and $K$ be
elements of $C$ such that $H \equiv_{\ge_{A, C}(\ga_0\jj \ga_1)} K$. Let us
write
$H$ in the form $H = \JJm{a_i \otimes b_i}{i < m}$ (the join in $C$, see
Proposition~\ref{P:sub-tensor}(iv)), with a positive integer $m$. For every
$i <
m$, there exists $a'_i \in A$ such that $\vv<a'_i, b_i> \in K$ and $a_i
\equiv_{\ga_0 \jj \ga_1} a'_i$. Since $a_i \equiv_{\ga_0 \jj \ga_1} a_i \mm
a_i'$
and $K$ is hereditary, we can replace $a'_i$ by $a_i \mm a_i'$; so we can
assume
that $a_i' \leq a_i$.  Thus there exist a positive integer $n$ and chains
 \[
   a_i  =  a_{i0} \geq a_{i1} \geq \cdots \geq a_{i, 2n} = a'_i
 \]
 such that for all $i < m$ and $j < n$, we have
 \begin{align*}
   a_{i, 2j} &\equiv_{\ga_0}a_{i, 2j + 1},\\
  a_{i, 2j + 1} &\equiv_{\ga_1}a_{i, 2j + 2}.
 \end{align*}
 Now, for all $j \leq  2n$, put $K_j = \JJm{a_{ij}\otimes b_i}{i < m}$ (the
joins
are formed in $C$). Note that $K_0 = H$ and $K_{2n} \ci K$. Furthermore,
for all $i
< m$ and $j < n$, we have
 \begin{align*}
   a_{i, 2j} \otimes b_i &\equiv_{\ge_{A, C}(\ga_0)} a_{i, 2j +
     1}\otimes b_i,\\
   a_{i, 2j + 1} \otimes b_i &\equiv_{\ge_{A, C}(\ga_1)}a_{i, 2j +
    2}\otimes b_i,
 \end{align*}
 from which it follows that $K_{2j} \equiv_{\ge_{A, C}(\ga_0)} K_{2j + 1}$ and
$K_{2j + 1} \equiv_{\ge_{A, C}(\ga_1)}K_{2j + 2}$. This proves the first
half of
the definition of $H \equiv_{\ge_{A, C}(\ga_0) \jj \ge_{A, C}(\ga_1)} K$;
the proof
of the other half is similar.
 \end{proof}

 \section{Tensor product of lattice congruences}\label{S:TensProd}

 Let $A$ and $B$ be lattices with zero and let $C$ be a sub-tensor product
of $A$
and~$B$.

 For $\ga \in \Con A$ and $\gb \in \Con B$, we define the \emph{$C$-tensor
product}
of $\ga$ and $\gb$ in~$C$, $\ga \odc{C} \gb$, by the formula
 \begin{align*}
   \ga \odc{C} \gb &= \ge_{A, C}(\ga) \mm \ge_{B, C}(\gb)\\
              &= (\ga \tim{C} \go_B) \mm (\go_A \tim{C} \gb).
 \end{align*}

 We write $\ga \odc{C} \gb$ in order to distinguish this congruence of $C$
from the
pure tensor $\ga \otimes \gb$ in the \jz-semilattice $\Con A \otimes \Con B$.

 In this section, we prove that $\ga \otimes \gb \mapsto \ga \odc{C} \gb$
extends
to a \jz-homomorphism $\ge_C \colon\Conc A \otimes \Conc B \to \Conc C$.  As a
first step, we prove that $\ga \otimes \gb \mapsto \ga \odc{C} \gb$ extends
to a
\jz-homomorphism $\ge_C \colon\Conc A \otimes \Conc B \to \Con C$. By
Proposition~\ref{P:TensUniv}, it is sufficient to prove the following:

 \begin{proposition}\label{P:OdotBimor} The map $\vv<\ga,\gb> \mapsto \ga
\odc{C}
\gb$ is a bimorphism from $\Con A \times \Con B$ to $\Con C$.
 \end{proposition}

 \begin{proof}
 Indeed, by Lemma~\ref{L:e0e1hom}, $\ge_{B, C}(\gw_B) = \go_C$, and so $\ga
\odc{C}
\go_B = \ge_{A, C}(\ga) \mm \ge_{B, C}(\go_B) = \ge_{A, C}(\ga) \mm \go_C =
\go_C$
and, similarly, $\go_A \odc{C} \gb = \go_{C}$.  Now compute:
 \begin{align*}
 \ga \odc{C} (\gb_0 \jj \gb_1)& =
 \ge_{A, C}(\ga) \mm \ge_{B, C} (\gb_0 \jj \gb_1)\\
 \intertext{(by Lemma~\ref{L:e0e1hom})}
 &= \ge_{A, C}(\ga) \mm (\ge_{B, C} (\gb_0) \jj \ge_{B, C}(\gb_1))\\
 \intertext{(by the distributivity of $\Con C$)}
 &= (\ge_{A, C}(\ga) \mm \ge_{B, C} (\gb_0)) \jj
    (\ge_{A, C}(\ga) \mm \ge_{B, C} (\gb_1))\\
 &= (\ga \odc{C} \gb_0) \jj (\ga \odc{C} \gb_1).
 \end{align*}

\vspace{-20pt} \end{proof}

 The crucial step in proving that $\ge_C$ maps $\Conc A \otimes \Conc B$ into
$\Conc C$ is the formula of Lemma~\ref{L:TensPpal}; we prepare its proof
with the
following statement:

 \begin{lemma}\label{L:CongInt}
 Let $\gg \in \Con C$ and let $b \leq  b'$ in $B$. Consider the following
subset
$\ga$ of $A$:
 \[
   \ga  =  \bigl \{\, \vv<x,y> \in A^2 \mid ((x \jj y)
    \otimes b) \jj((x \mm  y) \otimes b') \equiv_\gg (x \jj y)
     \otimes b'\,\bigr\}.
 \]
 Then $\ga$ is a congruence of $A$. \end{lemma}

 Note that $((x \jj y) \otimes b) \jj((x \mm  y) \otimes b')$, $(x \jj y)
\otimes
b' \in C$, so the formula makes sense.  We call this congruence $\ga$ the
\emph{$\vv<b, b'>$-projection of $\gg$ to $A$}.

 \begin{proof}
 We shall prove that $\ga$ satisfies the conditions listed in Lemma~I.3.8 of
\cite{Grat78}. It is obvious that $\ga$ is reflexive (because $x \otimes b'
= (x
\otimes b) \jj (x \otimes b')$ follows from $b \leq b'$) and $x \equiv_\ga
y$ iff
$x \mm  y \equiv_\ga x \jj y$, for all $x$, $y \in A$ (by the definition of
$\ga$).

Now let $x  \leq  y  \leq  z$ in $A$ such that $x \equiv_\ga y$ and $y
\equiv_\ga
z$. Then \begin{alignat*}{3}
    z \otimes b'& \equiv_\gg\, &&(y \otimes b') \jj_C (z \otimes b)
       &&\text{(since $y \leq z$ and $y \equiv_\ga z$)}\\ &\equiv_\gg&&(x
\otimes
b') \jj_C (y \otimes b) \jj_C (z \otimes b) \qq&&\text{(since $x \leq y$ and $x
\equiv_\ga y$)}\\
    &= &&(x \otimes b') \jj_C (z \otimes b)\qq &&
    \text{(since $y\leq z$)},
 \end{alignat*}
 thus $x \equiv_\ga z$.

 Finally, let $x$, $y$, $z \in A$ be such that $x  \leq  y$ and $x
\equiv_\ga y$.
We prove that $x \jj z \equiv_\ga y \jj z$ and $x \mm  z \equiv_\ga y \mm
z$. The
easier computation is for the join:
 \begin{alignat*}{3}
  (y \jj z) \otimes b'&  = \,&& (y \otimes b') \jj_C (z \otimes b')
       && \text{(by Proposition~\ref{P:OdotBimor})}\\
   &  = \,&& (y \otimes b') \jj_C (z \otimes b') \jj_C (z \otimes
   b) &&\text{(since $b \leq b'$)}\\
   &\equiv_\gg \,&& (x \otimes b') \jj_C (y \otimes b)
      \jj_C (z \otimes b') \jj_C (z \otimes b)\q
      &&\text{(since $x \equiv_\ga y$, $x \leq y$)}\\
    &  = \,&& ((x \jj z) \otimes b') \jj_C ((y \jj z) \otimes b)
     &&\text{(by Proposition~\ref{P:OdotBimor})},
 \end{alignat*}
 so that $x \jj z\equiv_\ga y \jj z$.

Now we compute the meet. Put
 \begin{align*}
   u  &=  ((x \mm  z) \otimes b') \jj_C ((y \mm  z) \otimes b),\\
   v  &=  (y \mm  z) \otimes b',\\
   u' &=  (x \otimes b') \jj_C (y \otimes b),\\
   v' &=  y \otimes b'.
 \end{align*}
 Note that $u  \leq  v$ and $u'  \leq  v'$. By definition, $x \equiv_\ga y$
means
that $u' \equiv_\gg v'$.  To prove that $x \mm  z \equiv_\ga y \mm  z$, we
have to
verify that $u \equiv_\gg v$; so it is sufficient to prove that the
interval $[u,
v]$ weakly projects up into the interval $[u', v']$, that is, $v
 \leq  v'$ and $u  =  v \mm  u'$.  It is obvious that $v
 \leq  v'$. To prove  $u  =  v \mm  u'$, compute, using the distributive law of
Proposition~\ref{P:sub-tensor}(vi):
 \begin{align*}
   v \mm  u'& = ((y \mm  z) \otimes b') \mm ((x \otimes b') \jj_C
    (y \otimes b))\\
   & = ((x \mm  z) \otimes b') \jj_C ((y \mm z) \otimes b)\\
   & = u. \end{align*}

\vspace{-20pt} \end{proof}

Next we show that the $C$-tensor product of two principal congruences is
principal
again.

 \begin{lemma}\label{L:TensPpal}
 Let $a  \leq  a'$ in $A$ and let $b  \leq  b'$ in $B$. Then $(a \otimes
b') \jj
(a' \otimes b)$, $a' \otimes b' \in C$ and the following formula holds:
 \[
   \gQ_A(a,a') \odc{C} \gQ_B(b,b')  =
     \gQ_{C}((a \otimes b') \jj (a' \otimes b), a' \otimes b').
 \]
 \end{lemma}

 \begin{proof}
 Put
 \begin{align*}
   \ga &= \gQ_A(a,a'),\\
   \gb &= \gQ_B(b,b'),\\
   \gg &= \gQ_{C}((a \otimes b') \jj_C (a' \otimes b),a' \otimes
       b').
 \end{align*}
 Note that
 \begin{align*}
   a \otimes b'
			   &\equiv_{\ge_{A, C}(\ga)} a' \otimes b',\\
   a' \otimes b
      &\equiv_{\ge_{B, C}(\gb)} a' \otimes b'.
 \end{align*}
 It follows that
 \[
    (a \otimes b') \jj (a' \otimes b) \equiv_{\ga \odc{C} \gb} a'
       \otimes b',
 \]
 whence $\gg \leq \ga \odc{C} \gb$.

Now let us prove the converse. The following statement will be helpful:

 \setcounter{claim}{0}
 \begin{claim}\label{C:1tens}
 Let $x  \leq  x'$ in $A$ and $y  \leq  y'$ in $B$ be such that
 $x \equiv_\ga x'$ and $y \equiv_\gb y'$. Then
 \[
   (x' \otimes y) \jj(x \otimes y') \equiv_\gg x' \otimes y'.
 \]
 \end{claim}

 \begin{proof}
 Let $\ga'$ be the $\vv<b, b'>$-projection of $\gg$ to $A$, see
Lemma~\ref{L:CongInt}. Then $a \equiv_{\ga'}a'$, by the definition
of~$\gg$, whence
$ \ga \leq \ga'$. Now let $u$, $u' \in A$ satisfy $u \leq u'$ and $u
\equiv_{\ga'}
u'$; define $\gb'_{u, u'}$ as the $\vv<u, u'>$-projection of $\gg$ to $B$.  Let
$\gb' = \MMm{\gb'_{u, u'}}{u, u' \in A,\ u \leq u',\ u \equiv_{\ga'} u'}$. Then
$\gb'$ is the intersection of a family of congruences of $B$, thus a
congruence of
$B$. By the definition of $\gb'$, we have $b \equiv_{\gb'}b'$; whence $\gb \leq
\gb'$. So if $x  \leq  x'$ in $A$ and $y  \leq  y'$ in $B$ such that $x
\equiv_{\ga'} x'$ and $y \equiv_{\gb'} y'$, then
 \[
   x' \otimes y' \equiv_\gg (x' \otimes y) \jj(x \otimes y').
 \]
 Since $\ga \leq \ga'$ and $\gb \leq \gb'$, the conclusion of the claim follows
immediately.
 \end{proof}

Now let $H$ and $K$ be elements of $C$ such that $H \equiv_{\ga \odc{C}
\gb}K$. We
write $H$ as
 \[
   H = \JJ\nolimits_C (x_i \otimes y_i \mid i < n),
 \]
 where $n$ is a positive integer, and $x_i \in A$, $y_i \in B$. Since $H
\equiv_{\ge_{A, C}(\ga) \mm \ge_{B, C}(\gb)}K$, for all $i<n$, there exist
$x^*_i
\leq  x_i$ and $y^*_i  \leq  y_i$ such that $x^*_i \equiv_{\ga}x_i$, $y^*_i
\equiv_{\gb}y_i$ and both $\vv<x^*_i, y_i>$ and $\vv<x_i, y^*_i>$ belong to
$K$.
But by Claim~\ref{C:1tens}, it follows that
 \[
   x_i \otimes y_i \equiv_{\gg} (x^*_i \otimes y_i) \jj(x_i
 \otimes y^*_i).
 \]
 By taking the join of these congruences over $i < n$, we obtain that $H
\equiv_{\gg} K'$ for some $K' \ci K$; by symmetry, we obtain the proof of the
symmetric inclusion and congruence, so that $H \equiv_{\gg}K$. Therefore, $\ga
\odc{C} \gb \leq \gg$. \end{proof}

 From Lemma~\ref{L:TensPpal} and Proposition~\ref{P:OdotBimor}, one can
then deduce
immediately the following statement:

 \begin{corollary}\label{C:TensComp}
 Let $\ga \in \Con A$ and $\gb \in \Con B$. If $\ga$ and $\gb$ are compact,
then
$\ga \odc{C} \gb$ is compact.
 \end{corollary}

 At this point, we have arrived at the existence of a (unique)
\jz-homomor\-phism
$\ge_C$ from $\Conc A \otimes \Conc B$ to $\Conc C$ such that, for all $\ga \in
\Conc A$ and all $\gb \in \Conc B$, we have
 \[
   \ge_C(\ga \otimes \gb)  =  \ga \odc{C} \gb.
 \]

\begin{example}\label{E:NonCon}
 We are relating $\Conc(A \otimes B)$ with $\Conc A \otimes \Conc B$
because, in
general, the lattice $\Con(A \otimes B)$ is not isomorphic to $\Con A
\otimes \Con
B$, not even if $A$ and $B$ are locally finite. Take $A = B = \go$, the
chain of
all non-negative integers. Then $\Con A = \Con B \iso \Pow \go$ (the power
set of
$\go$).  Since $\Con(A \otimes B)$ is an algebraic lattice and $\Con A
\otimes \Con
B$ is the \jz-semilattice of compact elements of an algebraic lattice, the
isomorphism $\Con A \otimes \Con B \iso \Con(A \otimes B)$ would imply that
in $A
\otimes B$ every congruence is compact, which is clearly not the case.

 We can be more specific. The isomorphism $\Con A \otimes \Con B \iso \Con(A
\otimes B)$ would imply that $\Con A \otimes \Con B$ is a complete lattice,
thus
that $\Pow\go\otimes\Pow\go$ is a complete lattice. However, $\Pow\go \otimes
\Pow\go$ is isomorphic to the \jz-subsemilattice of $\Pow(\go \times \go)$
generated by all \emph{rectangles}, that is, the subsets of the form $X
\times Y$
where $X$ and $Y$ are subsets of $\go$; in particular, it is not complete,
because
the set of all singletons $\{\vv<n, n>\}$, for $n \in \go$, does not have a
least
upper bound.

 This example contradicts Theorem 3.18 of \cite{GLQu81}.
 \end{example}

 \section{The Embedding Theorem}\label{S:Embedding}
 In this section, we will prove, still under the assumption that both $A$
and $B$
are lattices with zero and $C$ is a sub-tensor product of $A$ and $B$, that
the map
$\ge_C$ obtained in Section~\ref{S:TensProd} is a \jz-semilattice
\emph{embedding}.

 Our first lemma expresses the $\tim{C}$ operation on congruences by the
$\odc{C}$
operation:

 \begin{lemma}\label{L:timfrodot}
 Let $\ga \in \Con A$ and $\gb \in \Con B$. Then
 \[
   \ga\tim{C}\gb  = (\ga \odc{C} \gi_B)\vee(\gi_A \odc{C} \gb).
 \]
 \end{lemma}

 \begin{proof}
 By Lemma~\ref{L:e0e1hom}, we have $\ga \odc{C} \gi_B  = \ge_{A, C}(\ga) \mm
\ge_{B, C}(\gi_B) = \ge_{A, C}(\ga) = \ga \tim{C}\go_B$ and $\gi_A \odc{C}
\gb  =
\ge_{A, C}(\gi_A) \mm \ge_{B, C}(\gb) = \ge_{B, C}(\gb) = \go_A \tim{C}
\gb$. It
follows that $(\ga \odc{C} \gi_B) \jj (\gi_A \odc{C} \gb) \leq \ga \tim{C}\gb$.

 Conversely, let $H$ and $K$ be elements of $C$ such that
$H\equiv_{\ga\tim{C}\gb}K$. There exists a decomposition of $H$ in the form
 \[
   H = \JJ\nolimits_C (a_i \otimes b_i \mid i < n)
 \]
 with $n$ a positive integer and $a_i \in A$, $b_i \in B$. For all $i<n$, there
exist $a^*_i \leq  a_i$ and $b^*_i \leq  b_i$ such that $a^*_i \equiv_{\ga}
a_i$,
$b^*_i \equiv_{\gb} b_i$, and $\vv<a^*_i, b^*_i> \in K$. Thus
 \[
    a_i \otimes b_i \equiv_{\ga \odc{C} \gi_B} a^*_i \otimes
       b_i \equiv_{\gi_A \odc{C} \gb} a^*_i \otimes b^*_i,
 \]
 from which it follows that
 \begin{align*}
    H = \JJ\nolimits_C (a_i \otimes b_i \mid i < n)
        &\equiv_{\ga  \odc{C} \gi_B}
  \JJ\nolimits_C (a^*_i \otimes b_i \mid i < n)\\
   &\equiv_{\gi_A \odc{C}  \gb} \JJ\nolimits_C (a^*_i \otimes
      b^*_i \mid i < n)\\ & \ci  K.
 \end{align*}
 By symmetry, $H \equiv_{(\ga \odc{C} \gi_B) \jj (\gi_A \odc{C} \gb)}K$.
 \end{proof}

 \begin{lemma}\label{L:FundCont}
 Let $\ga$, $\ga' \in \Con A$ and $\gb$, $\gb' \in\Con B$. Then
 \[
   \ga \odc{C} \gb \leq \ga' \tim{C} \gb' \text{\q iff\q }
    \ga \leq \ga'\text{ or } \gb \leq \gb'.
 \]
 \end{lemma}

 \begin{proof}
 Let us assume that $\ga \odc{C} \gb \leq \ga' \tim{C} \gb'$ and $\ga \nleq
\ga'$,
$\gb \nleq \gb'$.  Then there are $a_0 < a_1$ in $A$ and $b_0 < b_1$ in $B$
such
that $a_0 \equiv_{\ga}a_1$ but $a_0\not \equiv_{\ga'}a_1$, and $b_0
\equiv_{\gb}b_1$ but $b_0\not \equiv_{\gb'}b_1$. It follows from
Lemma~\ref{L:TensPpal} that
 \[
    a_1 \otimes b_1 \equiv_{\ga \odc{C}  \gb}
     (a_0 \otimes b_1) \jj (a_1 \otimes b_0),
 \]
 thus, by assumption,
 \begin{equation}\label{Eq:abotnab}
   a_1 \otimes b_1 \equiv_{\ga' \tim{C} \gb'}
    (a_0 \otimes b_1) \jj (a_1 \otimes b_0).
 \end{equation}
 By Proposition~\ref{P:sub-tensor}(v), $(a_0 \otimes b_1) \jj (a_1 \otimes
b_0)=
(a_0 \otimes b_1) \uu (a_1 \otimes b_0)$. Thus, applying (\ref{Eq:abotnab}),
yields, for example, elements $x\in A$ and $y \in B$ such that $a_1
\equiv_{\ga'}x$,
$b_1\equiv_{\gb'}y$ and $\vv<x, y> \in a_0 \otimes b_1$. But since $a_1\not
\equiv_{\ga'}a_0$ and $b_1\not \equiv_{\gb'}b_0$, we have $x\ne 0$ and
$y\ne 0$.
Therefore, $x\leq a_0$ and $y\leq b_1$; whence, $a_1\equiv_{\ga'}a_0$, a
contradiction.  The reverse implication is trivial.
 \end{proof}

 \begin{lemma}\label{L:TensDist}
 The tensor product of distributive semilattices with zero is a distributive
semilattice with zero.
 \end{lemma}

\begin{proof}
 Let $A$ and $B$ be distributive semilattices with zero and let $I$ and $J$ be
bi-ideals of $A \times B$.  Set $X_0 = I \uu J$, and, for $n > 0$, let
$X_n$ be the
set of all $\vv<x, y> \in A \times B$ such that $\vv<x, y>$ is the lateral
join of
two elements of $X_{n - 1}$. It is obvious that $X_{n-1}\ci X_n$.

 \begin{all}{Claim 1}
 $X_n$ is a hereditary set, for all $n\geq 0$.
 \end{all}

 \begin{proof}
 The statement is obvious for $n = 0$.  Let us assume that it is true for
$n - 1$.
Let $\vv<u, v> \leq \vv<x, y> \in X_n$. By definition, $\vv<x, y>$ is a lateral
join of two elements of $X_{n - 1}$, that is, $\vv<x, y> = \vv<x_0, y> \jj
\vv<x_1,
y>$, where $\vv<x_0, y>$, $\vv<x_1, y> \in X_{n - 1}$ (or symmetrically).
Therefore, $u \leq x_0 \jj x_1$ in~$A$ and $v \leq y$ in $B$. Since $A$ is
distributive, there are $x_0' \leq x_0$ and $x_1' \leq x_1$ in $A$ such
that $u =
x_0' \jj x_1'$. Since $X_{n - 1}$ is hereditary,  $\vv<x_0', y>$,
$\vv<x_1', y> \in
X_{n - 1}$ and also $\vv<x_0', v>$, $\vv<x_1', v> \in X_{n - 1}$.  This implies
that $\vv<u,v>=\vv<x_0'\jj x_1',v>=\vv<x_0',v>\jj\vv<x_1',v>\in X_n$, since
$u =
x_0' \jj x_1'$ and the join, $\vv<x_0', v>  \jj \vv<x_1', v>$, is a lateral
join.
 \end{proof}

 \begin{all}{Claim 2}
 $I \jj J = \UUm{X_n}{n < \go}$.
 \end{all}

 \begin{proof}
 Obvious, from Claim 1.
 \end{proof}

 Now to prove the lemma, let $I$, $J$, and $K$ be bi-ideals of $A \times
B$.  Since
$(I \mm K) \jj (J \mm K) \ci (I \jj J) \mm K$, to prove distributivity, it is
enough to verify the reverse inclusion. So let $\vv<x, y> \in (I \jj J) \mm
K$.
Then $\vv<x,y> \in I \jj J$, so by Claim~2, $\vv<x, y> \in X_n$, for some
$n \geq
0$. We now prove by induction on $n$ that $\vv<x, y> \in (I \mm K) \jj (J
\mm K)$.
If $n = 0$, then $\vv<x, y> \in I \uu J$; since $\vv<x, y> \in K$, we
obtain that
$\vv<x, y> \in (I \uu J) \ii K \ci (I \ii K) \uu (J \ii K) \ci (I \mm K)
\jj (J \mm
K)$. Let us assume that the statement is true for $n - 1$ and let $\vv<x,
y> \in
X_n$. So, $\vv<x, y>$ is a lateral join, that is, $\vv<x, y> = \vv<x_0, y> \jj
\vv<x_1, y>$, where $\vv<x_0, y>$, $\vv<x_1, y> \in X_{n - 1}$ (or
symmetrically).
By the induction hypothesis, $\vv<x_0, y>$, $\vv<x_1, y> \in (I \mm K) \jj
(J \mm
K)$. Since $(I \mm K) \jj (J \mm K)$ is a bi-ideal and $\vv<x, y> =
\vv<x_0,y> \jj
\vv<x_1, y>$ is a lateral join, we conclude that $\vv<x, y> \in (I \mm K)
\jj (J
\mm K)$.
 \end{proof}

 We can also derive this lemma from the theory developed in F.~Wehrung
\cite{Wehr96}. Let $A$ and $B$ be distributive semilattices with zero. Thus
they
are conical refinement monoids in the sense of \cite{Wehr96}. By Theorem 2.7 in
\cite{Wehr96}, the tensor product $A \otimes^{\mathrm{cm}}B$ of $A$ and $B$
in the
category of commutative monoids, monoid homomorphisms, and monoid
bimorphisms (as
defined in Section 1 of \cite{Wehr96}) is a conical refinement monoid. But
it is
trivial that $A \otimes^{\mathrm{cm}}B$ is also the tensor product of $A$
and $B$
as defined in the present paper---this amounts to verifying that $A
\otimes^{\mathrm{cm}}B$ is a semilattice. For semilattices, distributivity is
equivalent to the refinement property and the refinement property is preserved
under tensor products (see \cite{Wehr96}), so the lemma follows.

\begin{remark}
 Even for finite lattices $A$ and $B$, one cannot deduce Lemma~\ref{L:TensDist}
directly from Theorem 3.3 of \cite{Fras76}, because the tensor product
considered
in \cite{Fras76} is the tensor product of arbitrary join-semilattices (not
necessarily with zero), and the resulting tensor product is not isomorphic
to ours,
in general. However, this difficulty is easy to overcome: if $A
\otimes^{\mathrm{F}}B$ is Fraser's tensor product of $A$ and $B$, then it
is easy
to see that $A \otimes B$, as defined in this paper, is the quotient of $A
\otimes^{\mathrm{F}}B$ by the bi-ideal generated by all elements of the form $x
\otimes^{\mathrm{F}}0_B$ ($x \in A$) and $0_A \otimes^{\mathrm{F}}y$ ($y
\in B$);
thus we can conclude Lemma~\ref{L:TensDist} by Theorem 3.3 of \cite{Fras76} for
finite lattices.
 \end{remark}

 Now we can state the embedding result:

 \begin{theorem}[\tbf{Embedding Theorem}]\label{T:Embedding}
 Let $A$ and $B$ be lattices with zero, and let $C$ be a sub-tensor product
of $A$
and $B$. Then the natural \jz-homomorphism
 \[
   \ge_C \colon \Conc A \otimes \Conc B \to \Conc C
 \]
 is a \jz-embedding.
 \end{theorem}

\begin{proof}
 Let $\gg =  \JJm{\ga_i \otimes\gb_i}{i < m}$ and $\gd = \JJm{\ga'_j \otimes
\gb'_j}{j < n}$ be elements of $\Conc A \otimes \Conc B$ (with the $\ga_i$,
$\ga'_j$ in $\Conc A$ and the $\gb_i$, $\gb'_j$ in $\Conc B$). We prove that
$\ge_C(\gg) \leq \ge_C(\gd)$ implies that $\gg \leq \gd$, which implies that
$\ge_C$ is an embedding.

 The assumption means that
 \[
   \JJm{\ga_i \odc{C}  \gb_i}{i < m} \leq  \JJm{\ga'_j \odc{C}
     \gb'_j}{j < n}.
 \]
 Now recall that $\ga'_j \odc{C}  \gb'_j = \ge_{A, C}(\ga'_j) \mm \ge_{B,
C}(\gb'_j)$. Using the fact that $\Con C$ is a distributive lattice and
that both
$\ge_{A, C}$ and $\ge_{B, C}$ are \jz-homomorphisms, it is easy to see that
this is
equivalent to saying that, for all $i < m$ and all $X \ci  n$, we have
 \begin{align*}
   \ga_i \odc{C}  \gb_i &\leq
    \ge_{A, C}\Bigl(\JJm{\ga'_j}{j \in X}\Bigr) \jj
    \ge_{B, C}\Bigl(\JJm{\gb'_j}{j \in n - X}\Bigr)\\
   &= \JJm{\ga_j'}{j \in X} \tim{C} \JJm{\gb_j'}{j \in n - X}\quad
       \text{(by Lemma~\ref{L:timfrodot})}.
 \end{align*}

 Therefore, by Lemma~\ref{L:FundCont}, for all $i < m$ and all $X \ci n$,
 \[
   \ga_i \leq  \JJm{\ga'_j}{j \in X}\text{\q or\q }
     \gb_i \leq  \JJm{\gb'_j}{j \in n -  X}.
 \]
 By Lemma~\ref{L:TensDist}, $\Conc A \otimes\Conc B$ is a distributive
semilattice,
thus $\Conc A \ootimes \Conc B$ is a distributive lattice. Therefore,
computing in
this lattice yields that, for all $i < m$,
 \begin{align*}
   \ga_i \otimes \gb_i& \leq
   \MMm{\JJm{\ga'_j \otimes \gi_B}{j \in X} \jj
   \JJm{\gi_A \otimes \gb'_j}{j \in n -  X}}{X \ci  n}\\
   & = \bigvee\left(\,(\ga'_j \otimes \gi_B)\mm (\gi_A
      \otimes \gb'_j) \mid j < n\,\right)\\
   & = \JJm{\ga'_j \otimes \gb'_j}{j < n} = \gd,
 \end{align*}
 whence, $\gg \leq \gd$.
 \end{proof}

 \section{The Isomorphism Theorem}\label{S:capped}
 We introduced capped bi-ideals in Definition~\ref{D:capping}.  We now
apply this concept to sub-tensor products and tensor products.

 \begin{definition}\label{D:SmTensLatt}
 Let $A$ and $B$ be lattices with zero. A \emph{capped sub-tensor product}
of $A$
and $B$ is a sub-tensor product $C$ of $A$ and $B$ such that every element
of $C$
is capped (that is, it is a finite \emph{union} of pure tensors).

 We say that the tensor product $A \otimes B$ is \emph{capped}, if every
compact
bi-ideal of $A \times B$ is capped.
 \end{definition}

 In this section, we will prove that if $C$ is a capped sub-tensor product
of $A$
and~$B$, then the embedding $\ge_C$ of the Embedding Theorem is an isomorphism.

 \begin{lemma}\label{L:SmoLat}
 Let $A$ and $B$ be lattices with zero. Then $A \otimes B$ is capped if and
only if
$A \otimes B$ is a capped sub-tensor product of $A \otimes B$. In
particular, if
every element of $A \otimes B$ is capped, then $A \otimes B$ is a lattice.
 \end{lemma}

\begin{proof}
 We prove the nontrivial direction. So, suppose that every element of $A
\otimes B$
is capped.  Let $H$ and $K$ be elements of $A \otimes B$. Then we can write
$H$ and
$K$ as
 \begin{align*}
   H &=  \UUm{a_i \otimes b_i}{i < m},\\
     K &=  \UUm{a'_j \otimes b'_j}{j < n},
 \end{align*}
 where $m$ and $n$ are positive integers and, for all $i < m$ and $j < n$,
$a_i$,
$a'_j \in A$ and $b_i$, $b'_j \in B$. Therefore, by Lemma~\ref{L:IntersTens},
 \[
   H \ii K =  \UUm{(a_i\mm  a'_j) \otimes(b_i\mm  b'_j)}
     {i < m\text{ and }j < n},
 \]
 whence, $H\ii K \in A \otimes B$. Thus $A\otimes B$ is a lattice. We
conclude the
argument by Proposition~\ref{P:AotBTensLatt}.
 \end{proof}

 We prepare the proof of the Isomorphism Theorem with the following statement:

 \begin{lemma}\label{L:CanoCon}
 Let $A$ and $B$ be lattices with zero, let $C$ be a capped sub-tensor
product of
$A$ and $B$. Let $n$ be a natural number, let $a \leq a'$ in $A$, $b \leq
b'$ in
$B$, and let $a_i \leq a'$ and $b_i \leq  b'$ (for all $i < n$) be such
that the
following element
 \[
   K = (a \otimes b') \uu (a' \otimes b) \uu \UUm{a_i \otimes
     b_i}{i < n}
 \]
 belongs to $C$.  Then the congruence $\gQ_C(K, a' \otimes b')$ is in the range
of~$\ge_C$.
 \end{lemma}

 \begin{proof}
 We prove this statement by induction on $n$. For $n = 0$, the congruence
$\gQ_C(K,
a' \otimes b')$ is, by Lemma~\ref{L:TensPpal}, equal to $\gQ_A(a, a') \odc{C}
\gQ_B(b, b')$, which belongs to the range of~$\ge_C$.

 Now assume that $n > 0$ and that the lemma holds for all integers less
than $n$.
We prove that we can assume, without loss of generality, that $a \leq  a_i
\leq
a'$, for all $i < n$. Indeed, let
 \[
   K' = (a \otimes b') \uu (a' \otimes b) \uu \UUm{(a \jj  a_i)
     \otimes b_i}{i < n}.
 \]
 It is obvious that $K \ci  K'$. On the other hand, for all $i < n$, we have
 \[
   (a \otimes b') \jj (a_i \otimes b_i) \ce (a \otimes b_i)
    \jj (a_i \otimes b_i) = (a \jj  a_i) \otimes b_i,
 \]
 from which it follows that
 \[
    K' \ci (a \otimes b') \jj (a' \otimes b) \jj  \JJm{a_i
      \otimes b_i}{i < n} = K,
 \] so that $K = K'$.

Similarly, we can assume that $b \leq  b_i \leq  b'$, for all $i < n$.

Set $a^\dagger = \JJm{a_i}{i < n}$.

\begin{all}{Claim 1}
 The elements $a^\dagger \otimes b'$, $K \ii (a^\dagger \otimes b')$ belong
to $C$
and the congruence $\gQ_C (a^\dagger \otimes b', K \ii (a^\dagger \otimes b'))$
belongs to the range of $\ge_C$.
 \end{all}

 \begin{proof}
 $a^\dagger \otimes b'$, $K \ii (a^\dagger \otimes b') \in C$ since pure
tensors
belong to $C$, $K \in C$ by assumption, and $C$ is closed under set
intersection,
by definition.

 For all $i < n$, using Proposition~\ref{P:sub-tensor}(i) and the fact that $b
\leq  b_i$, compute: \begin{equation}\label{Eq:Kint}
 \begin{aligned}
   K \ii (a_i \otimes b')& = (a \otimes b') \uu (a_i \otimes b)
    \uu \UUm{(a_j\mm  a_i) \otimes b_j}{j < n}\\
   & = (a \otimes b') \uu (a_i \otimes b_i) \uu \UUm{(a_j \mm
     a_i) \otimes b_j} {j < n,\ j \ne i}.
 \end{aligned}
 \end{equation}

 Thus, by the induction hypothesis and (\ref{Eq:Kint}), $\gQ_C(a_i \otimes
b', K
\ii (a_i \otimes b'))$ belongs to the range of $\ge_C$.

 Furthermore, for all $x \in A$, $\gQ_C(K \ii (x \otimes b'), x \otimes
b')$ is the
least congruence $\gQ$ of $C$ such that $[x \otimes b']{\gQ} \leq
[K]{\gQ}$. Since
$a^\dagger \otimes b' = \JJm{a_i \otimes b'}{i < n}$, it follows easily that
 \[
   \gQ_C(K\ii(a^\dagger \otimes b'), a^\dagger \otimes b') =
     \bigvee\bigl(\,\gQ_C(K \ii(a_i \otimes b'), a_i \otimes b')
     \mid i < n\,\bigr).
 \]
 The conclusion of the claim follows.
 \end{proof}

\begin{all}{Claim 2}
 The elements $a' \otimes b'$, $K$, $a^\dagger \otimes b'$, $K \ii (a^\dagger
\otimes b')$, $(a^\dagger \otimes b') \jj (a' \otimes b)$, belong to $C$
and the
following equation holds:
 \[
   \gQ_C(K, a' \otimes b') = \gQ_C(K \ii
    (a^\dagger \otimes b'), a^\dagger \otimes b') \jj \gQ_C((a^\dagger
\otimes b')
\jj (a'
    \otimes b), a' \otimes b').
 \]
 \end{all}

 \begin{proof}
 The elements listed belong to $C$ since pure tensors and mixed tensors
belong to
$C$, $K \in C$ by assumption, and $C$ is closed under set intersection.

 The inequality $\geq$ results immediately from the relations
 \[
   K \ci (a^\dagger \otimes b') \jj (a' \otimes b) \ci a'
    \otimes b'.
 \]
 Conversely, let $\gQ$ be the congruence on the right side of the equation.
Then
$(a^\dagger \otimes b') \jj  K \equiv_\gQ K$ and $(a^\dagger \otimes b')
\jj (a'
\otimes b) \equiv_\gQ a' \otimes b'$, thus we have
 \[
   a' \otimes b' = K \jj (a' \otimes b') \equiv_\gQ K \jj
     (a^\dagger \otimes b') \jj (a' \otimes b) \equiv_\gQ K \jj (a'
     \otimes b) = K,
 \]
 which proves the inequality $\leq$.
 \end{proof}

But it follows from Lemma~\ref{L:TensPpal} that
 \[
   \gQ_C((a^\dagger \otimes b') \jj (a' \otimes b), a' \otimes
        b') = \gQ_A(a^\dagger, a') \odc{C} \gQ_B(b, b'),
 \]
 thus this congruence belongs to the range of $\ge_C$.
 Therefore, it follows from Claims~1 and~2 that $\gQ_C(a' \otimes b', K)$
 belongs to the range of $\ge_C$.
\end{proof}

 \begin{theorem}[\tbf{Isomorphism Theorem}]\label{T:Isomorphism}
 Let $A$ and $B$ be lattices with zero, let $C$ be a capped sub-tensor
product of
$A$ and $B$. Then the natural embedding
 \[
   \ge_C \colon \Conc A \otimes \Conc B \to \Conc C
 \]
 is an isomorphism.

 In particular, if $C$ is a capped sub-tensor product of $A$ and $B$, then
 \[
   \Conc A\otimes\Conc B\iso\Conc C.
 \]
 \end{theorem}

 \begin{proof}
 By the Embedding Theorem and Lemma~\ref{L:SmoLat}, $\ge_C$ is a \jz-embedding;
it~remains to prove that $\ge_C$ is surjective. So, to conclude, it suffices to
prove that every $\gQ_C(H,K)$ (where $H$, $K \in C$) belongs to the range of
$\ge_C$. Without loss of generality, $H \ci  K$. Moreover, if $K = \JJm{a_i
\otimes
b_i}{i < n}$, then we have
 \[
   \gQ_C(H, K) =
    \JJm{\gQ_C((a_i \otimes b_i)\mm H, a_i \otimes b_i)}{i < n},
 \] so that it suffices to conclude in the case where $H \leq  K = a'
\otimes b'$,
for $a' \in A$ and $b' \in B$.

 Moreover, since $C$ is a capped sub-tensor product of $A$ and $B$, $H$ is
then a
finite \emph{union} of pure tensors $a_i \otimes b_i$ with $a_i \leq  a'$
and $b_i
\leq b'$. Hence, the theorem follows immediately from Lemma~\ref{L:CanoCon}
(with
$a = 0$ and $b = 0$).
 \end{proof}

 \section{Discussion}\label{S:Discussion}

 \subsection{Some corollaries} In this section, we list some consequences
of the
results of the last few sections.

 The following corollary is the first part of the Main Theorem of this paper as
stated in the Introduction:

 \begin{corollary}\label{C:mainfirstpart}
 Let $A$ and $B$ be lattices with zero.  If $A \otimes B$ is a lattice,
then there
is a natural embedding of $\Conc A \otimes \Conc B$ into $\Conc(A \otimes B)$.
 \end{corollary}
 \begin{proof}
 If $A \otimes B$ is a lattice, then $C = A \otimes B$ is a sub-tensor
product of
$A$ and $B$ by Proposition~\ref{P:AotBTensLatt}, so this corollary follows
from the
Embedding Theorem (Theorem~\ref{T:Embedding}).
 \end{proof}

 The next corollary is second part of the Main Theorem of this paper:

 \begin{corollary}\label{C:mainsecondpart}
 Let $A$ and $B$ be lattices with zero.  If $A \otimes B$ is capped, then
 \[
   \Conc A \otimes \Conc B \iso \Conc(A \otimes B).
 \]
 In fact, an isomorphism is exhibited by the natural map $\ge = \ge_{A
\otimes B}$.
 \end{corollary}

 \begin{proof}
 Indeed, if $A \otimes B$ is a capped tensor product, then $C = A \otimes
B$ is a
capped sub-tensor product of $A$ and $B$ by Lemma \ref{L:SmoLat}, so this
corollary
follows from the Isomorphism Theorem (Theorem~\ref{T:Isomorphism}).
 \end{proof}

 \begin{corollary}
 Let $A$ be a lattice with zero and let $S$ be a \emph{simple} lattice with
zero.
If $A \otimes S$ is capped, then $\Con A \iso \Con(A \otimes S)$.
 \end{corollary}

 \begin{proof}
 This is obvious because if $S$ is simple, then $\Con S = \Conc S$ is the
two-element chain and so $\Conc A \otimes \Conc S \iso \Conc A$.  By
Corollary~\ref{C:mainsecondpart}, $\Conc (A \otimes S) \iso \Conc A$, and
therefore, $\Con (A \otimes S) \iso \Con A$.
 \end{proof}

 We can also recover (and generalize) Theorem~4.2 and Corollary~4.4 of
\cite{GLQu81}, using the following trivial statement:

 \begin{proposition}\label{P:atom}
 Let $A$ and $B$ be \jz-semilattices. Then the atoms of $A\otimes B$ are
exactly
the pure tensors $a \otimes b$, where $a$ and $b$ are atoms of $A$ and $B$,
respectively.
 \end{proposition}

 \begin{corollary}
 Let $A$ and $B$ be lattices with zero with $|A|$, $|B| > 1$, let $C$ be a
capped
sub-tensor product of $A$ and $B$. Then $C$ is simple (resp., subdirectly
irreducible) if and only if $A$ and $B$ are simple (resp., subdirectly
irreducible).
 \end{corollary}

 \begin{proof}
 If $A$ and $B$ are simple, then, by the Isomorphism Theorem, $\Conc C \iso
\Conc A \otimes \Conc B$, and so $\Conc C$ is the two-element chain; it
follows that $C$ is simple.

 If $A$ and $B$ are subdirectly irreducible, then $A$ has a congruence
$\gF > \go$
with the property that $\gF \leq \ga$, for any congruence $\ga > \go$ of
$A$, and
$B$ has a congruence $\gQ > \go$ with the property that $\gQ \leq \gb$, for any
congruence $\gb > \go$ of $B$.  It is evident that $\gF \in \Conc A$ and
$\gQ \in
\Conc B$ and so $\gF \otimes \gQ$ is the unique atom of
$\Conc A \otimes \Conc B$
contained in all nonzero elements.  By the Isomorphism Theorem,
$\gF \odc{C} \gQ$ is
the unique atom of $\Conc C$ contained in all nonzero elements.  Since this
property is preserved when forming the ideal lattice, $\gF \odc{C} \gQ$ is
the unique atom of $\Con C$ contained in all nonzero elements, hence, $C$ is
subdirectly irreducible.
 \end{proof}

 \subsection{The paper \cite{GLQu81}}\label{S:GLQu81}
 The first draft of \cite{GLQu81} contained only the result stated in the
Introduction as the ``Main result of \cite{GLQu81}''.  The published version,
however, contained two generalizations:

 \begin{all}{Theorem 3.16 of \cite{GLQu81}}
 Let $A$ a finite lattice and let $B$ be an $A$-lower bounded lattice with $0$.
Then the isomorphism
 \[
   \Con A \otimes \Con B \iso \Con(A \otimes B)
 \]
 holds.
 \end{all}

 In this result, the following concept is used:

 \begin{definition}\label{D:lowerbounded}
 Let $A$ be a finite lattice and let $B$ be a lattice with $0$. We say that $B$
is \emph{$A$-lower bounded}, if, for every $n > 0$ and $n$-ary polynomial
$p_0$, any subset of $B$ of the form
 \[
   \setm{(p^{\mathrm{d}})_B(b_0, b_1, \dots, b_{n-1})}{p_A(a_0,
     a_1, \dots, a_{n-1}) = (p_0)_A(a_0, a_1, \dots, a_{n-1})}
 \]
 (where $p$ ranges over all $n$-ary polynomials) has a largest element.
 \end{definition}

 By Lemma~\ref{L:IntersTens}(iv), we have a formula for the meet of two
elements of the tensor product of lattices.  Unfortunately, the right side is,
in general, an infinite union.  However, if $B$ is $A$-lower bounded, then, for
given $p(\vec a)$ and $q(\vec c)$, we can choose the largest
$p^{\mathrm{d}}(\vec b)$ and $q^{\mathrm{d}}(\vec d)$ and so the right side
equals a finite subunion.

 Therefore, the condition that $A$ be finite and $B$ be $A$-lower bounded
is the
most natural one under which $A \otimes B$ is a lattice.  Thus Theorem 3.16 of
\cite{GLQu81} is a natural extension of the Main result of \cite{GLQu81}.
Unfortunately, the proof retained from the finite case the assumption that $B$
has a $1$, so, in fact, the result is proved in \cite{GLQu81} only under the
additional assumption that $B$ have a unit element.

 The above discussion shows that if $A$ is finite and $B$ is $A$-lower bounded,
then $A \otimes B$ is capped, so our Isomorphism Theorem proves this result
without any additional assumptions.

 Then \cite{GLQu81} goes on to argue that nothing changes if $A$ is only
assumed to be locally finite:

 \begin{all}{Theorem 3.18 of \cite{GLQu81}}
 Let $A$ be a locally finite lattice with zero and let $B$ be an $A$-lower
bounded lattice with $0$. Then the isomorphism
 \[
   \Con A \otimes \Con B \iso \Con(A \otimes B)
 \]
 holds.
 \end{all}

 This is obviously not true.  The proof of Theorem~3.16 computes the compact
elements of $\Con(A \otimes B)$.  If $A$ is only locally finite, then these
computations show very little about congruences in general, see our
Example~\ref{E:NonCon} (in the example, both $A$ and $B$ are distributive, so
$B$ is trivially $A$-lower bounded).  The correct form of Theorem~3.18 of
\cite{GLQu81} switches to the isomorphism:
 \[
   \Conc A \otimes \Conc B \iso \Conc(A \otimes B),
 \]
 which indeed follows from our Isomorphism Theorem.

 \subsection{The papers \cite{Farl96} and \cite{GrSc94}}\label{S:GSFa}
 Let $L$ be a lattice and let $D$ be a bounded distributive lattice.  The
lattice $L[D]$ is defined in G. Gr\"atzer and E. T. Schmidt \cite{GrSc94} as
follows. First, if $D$ is finite, then let $P$ be the poset of join-irreducible
elements of $D$ and define $L[D]$ as the \emph{function lattice} $L^P$,
that is,
the lattice of all order-preserving maps from $P$ to $L$, partially ordered
componentwise. For an arbitrary bounded distributive lattice $D$, define $L[D]$
as the \emph{direct limit} of all $L[D']$, where $D'$ is a finite
$\set{0,1}$-sublattice of $D$, with the natural embeddings. This construction
yields a lattice $L[D]$ that is isomorphic to the lattice studied in
\cite{GrSc94}, see Lemma~1 in \cite{GrSc94}. This construction is studied from
a more topological point of view in \cite{Farl96}.

 $L \otimes D$ is defined, if $L$ has a zero; one would expect that, in this
case, the isomorphism $L[D] \iso L \otimes D$ holds. The reality is slightly
more awkward. For a lattice $K$, let $K^\dd$ denote the \emph{dual lattice}
of $K$.  Then the following isomorphism holds:
 \begin{equation}\label{Eq:duals}
 L[D]\iso(L^\dd\otimes D)^\dd,
 \end{equation}
provided that $L$ has a \emph{greatest} (as opposed to least) element.
Formula~(\ref{Eq:duals}) is easy to establish for finite $D$ and the general
case follows by a direct limit argument.

In both papers, \cite{Farl96} and \cite{GrSc94}, the congruence lattice of
$L[D]$ is computed from $\Con L$ and $\Con D$. Note that $\Con D$ is isomorphic
to the ideal lattice of the generalized Boolean algebra generated by $D$. For
example, Theorem~4 of \cite{GrSc94} states that
 \begin{equation}\label{Eq:ConcLD}
 \Conc L[D]\iso(\Conc L)[\Conc D].
 \end{equation}
If $L$ is a lattice with a greatest element, then the isomorphism in
(\ref{Eq:ConcLD}) is an elementary consequence of the Isomorphism Theorem,
although the proof requires a number of tedious translations between the two
constructions.

We cannot directly deduce Formula (\ref{Eq:ConcLD}) for an arbitrary lattice
$L$ from the Isomorphism Theorem.  However, this is possible. In \cite{GrWe2},
there is even a generalization of Formula (\ref{Eq:ConcLD}) for an arbitrary
lattice $L$  and for an arbitrary (not necessarily distributive) bounded
lattice $D$, with an analogue of Formula (\ref{Eq:ConcLD}).

 \subsection{Some open problems}
 Many of the problems asking whether the conditions we use in this paper
are also necessary are still open.

 If $A$ is a locally finite lattice and $B$ is $A$-lower bounded, then it
is easy
to see that $A \otimes B$ is capped.

 \begin{problem}
 Do there exist lattices $A$ and $B$ with zero so that $A \otimes B$ is
capped and
neither $A$ nor $B$ is locally finite?
 \end{problem}

 \begin{problem}\label{P:SmTensLatt}
 Let $A$ and $B$ be lattices with zero, let $C$ be a sub-tensor product of
$A$ and
$B$. Is $C$ a capped sub-tensor product of $A$ and $B$?
 \end{problem}

 Specializing Problem~\ref{P:SmTensLatt} to $C = A \otimes B$ yields the
following
question:

 \begin{problem}\label{P:SmPair1}
 Let $A$ and $B$ be lattices with zero. If $A \otimes B$ is a lattice, is $A
\otimes B$ capped?
 \end{problem}

 We prove in \cite{GrWe1} that the answer to Problem~\ref{P:SmPair1} is
positive,
if $A$ or $B$ is locally finite.  Moreover, we provide the example $M_3 \otimes
\tup{F}(3)$ of a tensor product of lattices with zero that is not a
lattice.  In
\cite{GrWe4} we exhibit a three-generated planar lattice $L$ such that $M_3
\otimes
L$ is not a lattice.

 If $A \otimes B$ is capped, then $C = A \otimes B$ is the largest sub-tensor
product of $A$ and $B$.

 \begin{problem}\label{P:SmPair2}
 Let $A$ and $B$ be lattices with zero. Does there always exist a largest
sub-tensor product (resp., a largest capped sub-tensor product) of $A$ and $B$?
 \end{problem}

 For a join-semilattice $S$ with zero, denote by $\ConL S$ the set of all
L-congru\-ences of $S$.  The results of Section~\ref{S:L-congruences}
suggest that
$\vv<\ConL S, \ci>$ must behave to some extent as the congruence lattice of a
lattice.

 \begin{problem}
 What is the structure of $\vv<\ConL S, \ci>$?
 \end{problem}

 Theorem~4.6 of \cite{GLQu81} investigates subdirectly irreducible
quotients of a
tensor product of lattices with zero. For finite lattices $A$ and $B$, it
is proved
that the completely meet-irreducible congruences of $A \otimes B$ are
exactly the
congruences of the form $\ga \congtimes \gb$, where $\ga$, $\gb$ are completely
meet-irreducible congruences of $A$, $B$, respectively. The proof of
Theorem~4.6
(and consequently, of Theorem~4.7 and Theorem~4.9 of \cite{GLQu81}) does
not apply
to infinite lattices.

 \begin{problem}
 Let $A$ and $B$ be lattices with zero, let $C$ be a capped sub-tensor
product of
$A$ and $B$. Are the completely meet-irreducible congruences of $A \otimes B$
exactly the congruences of the form $\ga \tim{C} \gb$ where $\ga$, $\gb$ are
completely meet-irreducible congruences of $A$, $B$, respectively?
 \end{problem}

 Note that by Corollary~\ref{C:NabCon}(iii) and the Isomorphism Theorem, every
congruence of the form $\ga \tim{C} \gb$, where $\ga$, $\gb$ are completely
meet-irreducible congruences of $A$, $B$, respectively, is a
completely meet-irreducible congruence of $C$.

 \begin{notation}
 For a lattice $L$ with zero,
 \begin{enumerate}
 \item let $\tup{D}(L)$ denote the maximal distributive quotient of $L$;
 \item let $\gQ_{\tup{D}}(L)$ denote the kernel of the natural homomorphism
onto
the maximal distributive quotient of $L$;
 \item let $\gQ_{\tup{M}}(L)$ denote the kernel of the natural homomorphism
onto
the maximal modular quotient of $L$.
 \end{enumerate}
 \end{notation}

 \begin{problem}
 Let $A$ and $B$ be lattices with zero. If $A \otimes B$ is capped, is it
then true
that
 \[
   \tup{D}(A \otimes B) \iso \tup{D}(A) \otimes \tup{D}(B)
 \]
 holds?
 \end{problem}

The proof presented in Theorem~4.7 of \cite{GLQu81} applies to the finite case.

\begin{problem}
 Let $A$ and $B$ be lattices with zero.  If $A \otimes B$ is capped, is it then
true that
 \[
  \gQ_{\tup{M}}(A \otimes B) = (\gQ_{\tup{D}}(A) \congtimes
  \gQ_{\tup{M}}(B)) \mm (\gQ_{\tup{M}}(A) \congtimes
  \gQ_{\tup{D}}(B))
 \]
 holds?
 \end{problem}
 Again, the proof presented in Theorem~4.7 of \cite{GLQu81} applies in the
finite
case.

 \section*{Acknowledgment}
 This work was partially completed while the second author was visiting the
University of Manitoba. The excellent conditions provided by the Mathematics
Department, and, in particular, a very lively seminar, were greatly
appreciated.

\end{document}